\documentclass[5p,twocolumn]{elsarticle}

\usepackage{lineno,hyperref}
\usepackage{stfloats}
\usepackage{tikz} 
\usepackage{amsmath}
\usepackage{mathtools}
\usepackage{amsfonts}
\usepackage{float}
\usepackage{lmodern,bm}                
\usepackage[T1]{sansmath} 
\usepackage{amsthm}
\usepackage{url}
\usepackage[toc,page]{appendix}
\usepackage{graphicx}
\usepackage{subcaption}
\usepackage{verbatim} 
\usetikzlibrary{shapes.multipart}
\usetikzlibrary{matrix}
\usetikzlibrary{positioning}
\usepackage{float}
\usepackage{color}
\usepackage{booktabs}
\usepackage{placeins}
\usepackage{soul}

\DeclarePairedDelimiter\floor{\lfloor}{\rfloor}

\journal{Journal of Quantitative Spectroscopy \& Radiative Transfer}

\begin{document}

\begin{frontmatter}

\title{Calder\'on preconditioning of PMCHWT boundary integral equations for scattering by multiple absorbing dielectric particles}

\author[ucl]{Antigoni Kleanthous}
\ead{antigoni.kleanthous.12@ucl.ac.uk}
\author[ucl]{Timo Betcke}
\author[ucl]{David P. Hewett}
\author[ucl]{Matthew W. Scroggs}
\author[met,hertfordshire]{Anthony J. Baran}

\address[ucl]{University College London, Department of Mathematics, London, UK}
\address[met]{Met Office, Exeter, UK}
\address[hertfordshire]{University of Hertfordshire, School of Physics, Astronomy, and Mathematics, Hertfordshire, UK}

\begin{abstract}
We consider the simulation of electromagnetic scattering by single and multiple isotropic homogeneous dielectric particles using boundary integral equations. Galerkin discretizations of the classical Poggio-Miller-Chang-Harrington-Wu-Tsai (PMCHWT) boundary integral equation formulation provide accurate solutions for complex particle geometries, but are well-known to lead to ill-conditioned linear systems. 
In this paper we carry out an experimental investigation into the performance of Calder\'on preconditioning techniques for single and multiple absorbing obstacles, which involve a squaring of the PMCHWT operator to produce a well-conditioned second-kind formulation. For single-particle scattering configurations we find that Calder\'on preconditioning is actually often outperformed by simple ``mass-matrix'' preconditioning, i.e.\ working with the strong form of the discretized PMCHWT operator. 
In the case of scattering by multiple particles we find that 
a significant saving in computational cost can be obtained by performing block-diagonal Calder\'on preconditioning in which only the self-interaction blocks are preconditioned. Using the boundary element software library Bempp (\url{www.bempp.com}) the numerical performance of the different methods is compared for a range of wavenumbers, particle geometries and complex refractive indices relevant to the scattering of light by atmospheric ice crystals. 
\end{abstract}



\end{frontmatter}

\section{Introduction}
Calculating the radiative effects of cirrus clouds is a central problem in climate modelling \cite{liou2016light,baran2012single}. Cirrus clouds appear at high altitudes (usually greater than 6km \cite{baran2012single}), and their effect on climate is significant, since at any given time they cover around 30\% of the Earth's surface at mid-latitudes, and 60-80\% in the tropics, and are not confined to a particular latitude or season \cite{baran2009review}. The ice crystals in cirrus clouds vary considerably in size and shape, and are generally highly non-spherical, taking forms such as hexagonal columns, hexagonal plates and bullet rosettes \cite{baran2012single}, and aggregates of these, and the accurate simulation of electromagnetic scattering by such ice crystals is therefore a challenging problem \cite{liou2016light,baran2012single}. 

A number of different approaches are available for the simulation of electromagnetic dielectric scattering, each with their own advantages and disadvantages. 
For particles of small to moderate size (relative to the wavelength) there are ``numerically exact'' methods \cite{mishchenko2014electromagnetic} such as the Discrete Dipole Approximation (DDA) \cite{draine1994discrete, yurkin2007discrete}, the Finite-Difference Time-Domain (FDTD) \cite{yang1996finite, yang2000finite, yurkin2007systematic} and Pseudo-Spectral Time-Domain (PSTD) \cite{liu2012comparison} methods, and the  Extended Boundary Condition \cite{havemann2001extension,kahnert2013t, baran2001calculation} and Invariant Imbedded  \cite{bi2013efficient, bi2014accurate} T-matrix methods. 
For particles of large size (relative to the wavelength) one can use ``approximate'' high-frequency methods such as Geometric Optics/ray tracing and the Kirchhoff approximation (see e.g.\ \cite{muinonen1989scattering, muinonen1996light,
macke1996single, yang1996geometric, mishchenko1998incorporation,
liou2002introduction, borovoi2003scattering, hesse2008modelling, bi2009simulation,
bi2011scattering, hesse2012modelling, bi2014assessment}). Scattering by multiple particles is also a well-studied problem. For a general overview of the classical multi-particle scattering literature we refer the reader to \cite{martin2006multiple}. We also mention the T-matrix based methods in \cite{koc1998calculation, gumerov2005computation, zhang2007fast} and the fast solvers presented recently in \cite{bremer2015high, ganesh2015efficient}.

The use of boundary integral equations and their discretization using boundary element methods (BEMs) is well-established in the electrical engineering community but has only recently started getting serious attention in the atmospheric physics community \cite{groth2015boundary,baran2017application,yu2017electromagnetic}. BEM is a ``numerically exact'' method which can provide highly accurate simulations of scattering by complex ice crystal shapes. It has no inherent restriction on the complexity of the scatterer, and is also well-suited to higher frequency problems, since it is based on a reformulation of the scattering problem as a boundary integral equation on the scatterer's surface, which reduces the dimensionality of the domain to be discretized from an unbounded 3D domain to a bounded 2D domain. Early applications of BEM to the simulation of light scattering by simple ice crystals include \cite{mano2000exact} and \cite{nakajima2009development}. More comprehensive studies of complex crystal shapes (including hexagonal columns with conventional and stepped cavities, bullet rosettes and Chebyshev ice particles) have been given recently by Groth \textit{et al.} \cite{groth2015boundary} and Baran and Groth \cite{baran2017application}, using the open-source BEM software library Bempp \cite{smigaj2015solving}, available at \url{www.bempp.com}. 
In \cite{groth2015boundary,baran2017application} it was shown (by comparison with a T-matrix method) that, using a discretization with at least 10 boundary elements per wavelength, Bempp can compute far-field quantities with 1\% relative error, with reciprocity also satisfied to a similar accuracy. The results in \cite{groth2015boundary,baran2017application} were obtained using a standard desktop machine and were limited to size parameters up to 15, but Bempp supports parallelization and hence with HPC architectures much larger problems can be tackled. 
We end this brief literature review of BEM in atmospheric physics applications by mentioning recent work by Yu \textit{et al.} \cite{yu2017electromagnetic}, in which boundary integral equations were applied to scattering by multiple dielectric particles illuminated by unpolarized high-order Bessel vortex beams, and also by Groth \textit{et al.} \cite{groth2018hybrid}, where a ``hybrid numerical-asymptotic'' BEM was presented for 2D high frequency scalar dielectric scattering problems which achieves fixed accuracy with frequency-independent computational cost.

A number of different boundary integral equation formulations are available for dielectric scattering problems, the most popular of which (and the one used in \cite{groth2015boundary}) is the PMCHWT formulation due to Poggio, Miller, Chang, Harrington, Wu and Tsai \cite{poggio1970integral, wu1977scattering, mautz1977electromagnetic,harrington1989boundary}. 
The PMCHWT formulation provides accurate solutions for complex scatterer geometries, but is well-known to suffer from ill-conditioning, requiring a large number of iterations when the associated linear system is solved using an iterative method such as the Generalized Minimal Residual Method (GMRES)\cite{saad1986gmres}. 
One therefore requires a preconditioning strategy. In \cite{groth2015boundary,baran2017application}, algebraic preconditioning was applied, based on $\mathcal{H}$-matrices. But this incurs a large memory overhead. An alternative operator-based approach, studied e.g.\ in 
\cite{contopanagos2002well, christiansen2002preconditioner, antoine2008integral, bagci2009calderon, cools2011calderon, yan2010comparative, niino2012calderon, boubendir2015integral}, is ``Calder\'on preconditioning'', which involves squaring the original boundary integral operator to regularize the problem at the continuous level, i.e.\ before discretization. This accelerates GMRES convergence, at the expense of an increased cost per GMRES iteration. 

In this paper, we investigate the performance of Calder\'on preconditioning strategies in the context of light scattering by single and multiple complex ice crystals. We also investigate another operator-based preconditioning strategy, ``mass-matrix preconditioning'', which uses the discrete strong form of the PMCHWT operator \cite{kirby2010functional,betcke2017product} and does not require a second application of the PMCHWT operator.
For the case of scattering by multiple particles we consider three preconditioning strategies: mass matrix preconditioning; full Calder\'on preconditioning, involving a second application of the PMCHWT operator; and diagonal Calder\'on preconditioning, in which only the self-interaction blocks are preconditioned. The numerical performance of the three methods is compared for a range of wavenumbers, particle geometries and complex refractive indices. 
Our numerical simulations are carried out in Bempp, and the computational cost of each preconditioning strategy is measured by the total number of matrix-vector products (each corresponding to an application of one boundary integral operator) required to solve the discretized system to a specified tolerance.

Our main findings are that: (i) for scattering by a single particle, mass-matrix preconditioning is more effective than Calder\'on preconditioning; and (ii) for scattering by multiple particles, block-diagonal Calder\'on preconditioning outperforms both mass-matrix preconditioning and full Calder\'on preconditioning.

The paper is organized as follows. In Section \ref{sct:scattering_problem}, we outline the problem of scattering by multiple  homogeneous, isotropic, dielectric particles. In Section \ref{sct:operators}, we review the key properties of the electromagnetic potential operators and boundary integral operators arising in the boundary integral equation formulation. In Section \ref{sct:PMCHWT}, we formulate the PMCHWT boundary integral equation in the setting of scattering by multiple particles. 
In Section \ref{sct:preconditioning}, we describe our Galerkin discretization and the different preconditioning strategies we study. In Section \ref{sct:complexity}, we compare the computational complexity of the different strategies. In Section \ref{sct:benchmarks}, we present benchmarking results for simple particle shapes, investigating the performance of the preconditioners in various configurations of single and multiple particles, as a function of key parameters such as the number, size, separation and material properties of the constituent particles.
In Section \ref{sct:ice_crystals}, we apply our solver to a range of scattering problems relevant to light scattering by atmospheric ice crystals. Example Python notebooks are available on \url{www.bempp.com}. Concluding remarks are given in Section \ref{sct:conclusion}.

\section{The scattering problem}\label{sct:scattering_problem}

\begin{figure*}[ht!]
\centering
\includegraphics[width = 0.9\textwidth]{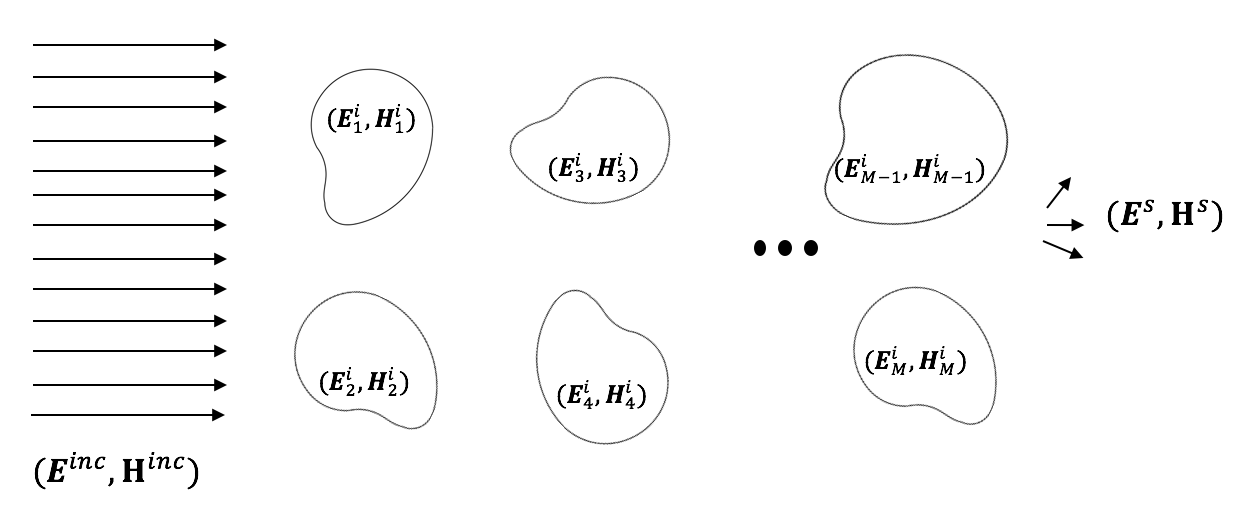}
\caption{Scattering by multiple particles.}
\label{fig:multiple_problem}
\end{figure*}

We consider electromagnetic scattering by a collection of $M$ disjoint arbitrary 3D isotropic homogeneous dielectric scatterers occupying bounded domains $\Omega^i_m\subset \mathbb{R}^3$, $m=1,\ldots,M$, with boundaries $\Gamma_m = \partial \Omega^i_m$, in a homogeneous exterior medium $\Omega^e = \mathbb{R}^3 \backslash \overline{\cup_{m=1}^M\Omega^i_m}$, as in Figure \ref{fig:multiple_problem}. 
The electric and magnetic fields in the interior domains $\Omega^i_m$, $m=1,\ldots,M$, and the exterior domain $\Omega^e$, will be denoted $(\mathbf{E}^i_m, \mathbf{H}^i_m)$ and $(\mathbf{E}^e, \mathbf{H}^e)$ respectively. They are assumed to satisfy the time-harmonic Maxwell equations
\begin{alignat}{3}
\label{eqn:maxwell_ei}\nabla \times \mathbf{E}^i_m &=& i \omega \mu_m \mathbf{H}^i_m, &\qquad\text{in }\Omega^i_m,\,m=1,\ldots,M, \\
\label{eqn:maxwell_Hi}\nabla \times \mathbf{H}^i_m &=& - i \omega \epsilon_m \mathbf{E}^i_m, &\qquad\text{in }\Omega^i_m,\,m=1,\ldots,M,
\end{alignat}
and
\begin{alignat}{3}
\label{eqn:maxwell_ee}\nabla \times \mathbf{E}^e &=& i \omega \mu_e \mathbf{H}^e, &\qquad\text{in }\Omega^e,\\
\label{eqn:maxwell_He}\nabla \times \mathbf{H}^e &=& - i \omega \epsilon_e \mathbf{E}^e,&\qquad\text{in }\Omega^e,
\end{alignat}
together with the transmission boundary conditions
\begin{alignat}{3}
\label{bc1}\mathbf{E}^i_m (\mathbf{x}) \times \mathbf{n} &= \mathbf{E}^e (\mathbf{x}) \times \mathbf{n}, && \quad \mathbf{x} \in \Gamma_m, \ m=1,\ldots,M,\\
\label{bc2}\mathbf{H}^i_m (\mathbf{x}) \times \mathbf{n} &= \mathbf{H}^e (\mathbf{x}) \times \mathbf{n}, &&\quad  \mathbf{x} \in \Gamma_m, \ m=1,\ldots,M.
\end{alignat}
Here we assume a time-dependence of the form $e^{-i\omega t}$, with angular frequency $\omega >0$. The parameters $\epsilon_m$, $\epsilon_e$ and $\mu_m$, $\mu_e$, represent respectively the electric permittivity and the magnetic permeability of the domains, and $\mathbf{n}$ is the unit normal vector on $\Gamma_m$ pointing into $\Omega^e$. 

In the scattering problem, an incident field $(\mathbf{E}^{inc},\mathbf{H}^{inc})$ (for instance, a plane wave) gives rise to internal fields $(\mathbf{E}^i_m,\mathbf{H}^i_m)$ in $\Omega^i_m$ and a scattered field $(\mathbf{E}^s,\mathbf{H}^s)$ in the exterior domain $\Omega^e$. The latter is assumed to satisfy the Silver-M\"{u}ller radiation condition, and the total exterior field is then the sum of incident and scattered fields
\begin{alignat}{3}
&\mathbf{E}^e &=& \mathbf{E}^{inc} &+ \mathbf{E}^s, \qquad&\text{in }\Omega^e,\\
&\mathbf{H}^e &=& \mathbf{H}^{inc} &+ \mathbf{H}^s, \qquad&\text{in }\Omega^e.
\end{alignat}

It is sufficient to solve for either the electric or magnetic fields and then recover the remaining fields by (\ref{eqn:maxwell_ei})-(\ref{eqn:maxwell_Hi}) and (\ref{eqn:maxwell_ee})-(\ref{eqn:maxwell_He}). In what follows, we will solve for the electric fields $\mathbf{E}^i_m$, $\mathbf{E}^e$, which satisfy
\begin{alignat}{3}
\label{single_equation_i}&\nabla \times (\nabla \times \mathbf{E}^i_m) &-& k_m^2 \mathbf{E}^i_m &= 0,\qquad&\text{in }\Omega^i_m,\\
\label{single_equation_e}&\nabla \times (\nabla \times \mathbf{E}^e) &-& k_e^2 \mathbf{E}^e &= 0,\qquad&\text{in }\Omega^e,
\end{alignat}
where $k_m = \omega \sqrt{\mu_m \epsilon_m}$ and $k_e = \omega \sqrt{\mu_e \epsilon_e}$ are the wavenumbers in the respective domains.

\section{Boundary integral operators}\label{sct:operators}
The basic idea behind the BEM is that the fields $(\mathbf{E}^i_m,\mathbf{E}^s)$ can be represented as potentials whose densities are the (unknown) boundary traces $(\gamma^-_{D,m}\mathbf{E}^i_m,\gamma^-_{N,m}\mathbf{E}^i_m, \gamma^+_{D,m}\mathbf{E}^s, \gamma^+_{N,m}\mathbf{E}^s)$ (see Eqns \eqref{eqn:StrattonChu_int}-\eqref{eqn:StrattonChu_ext} below). The transmission conditions imply that these traces satisfy certain boundary integral equations, which can be solved numerically using a BEM. 
In this section, we briefly recall the basic definitions of the potentials and boundary integral operators that underpin the BEM. For further details, including a discussion of the function space setting, the reader is directed to \cite{buffa2003galerkin}. 

Given a wavenumber $k$ and a bounded Lipschitz open set $\Omega$ with boundary $\Gamma=\partial\Omega$ and outward unit normal vector $\mathbf{n}$, we define the electric and magnetic potential operators, applied to a boundary vector field $\mathbf{v}$, by 
\begin{alignat}{2}
&\mathcal{E} \mathbf{v}(\mathbf{x}) &:=& \mathrm{i}k \int _{\Gamma} \mathbf{v}(\mathbf{y}) G (\mathbf{x}, \mathbf{y}) d \Gamma (\mathbf{y}) \nonumber \\
&&&- \frac{1}{\mathrm{i}k} \nabla_\mathbf{x} \int _\Gamma \nabla_\mathbf{y} \cdot \mathbf{v}(\mathbf{y}) G (\mathbf{x}, \mathbf{y}) d\Gamma (\mathbf{y}), \\
&\mathcal{H} \mathbf{v}(\mathbf{x}) &:=& \nabla_\mathbf{x} \times \int _\Gamma \mathbf{v}(\mathbf{y}) G (\mathbf{x}, \mathbf{y}) d \Gamma (\mathbf{y}),
\end{alignat}
where $G (x,y) = \displaystyle\frac{\exp (\mathrm{i}k | \mathbf{x} - \mathbf{y}|)}{4 \pi | \mathbf{x} - \mathbf{y}|}$.
On $\Gamma$ we define the interior ($-$) and exterior ($+$) Dirichlet (tangential) and Neumann traces $\gamma_{D}^{\pm}$, $\gamma_{N}^{\pm}$ (for details see e.g. \cite{buffa2003galerkin}), which for smooth vector fields $\mathbf{u}^+$ defined in $\mathbb{R}^3 \backslash \overline{\Omega}$, and $\mathbf{u}^-$ defined in $\Omega$, satisfy 
\begin{alignat}{3}
\gamma_{D}^{\pm} \mathbf{u}^\pm (\mathbf{x}) &= \mathbf{u}^\pm(\mathbf{x}) \times \mathbf{n}(\mathbf{x}), \quad &\mathbf{x} \in \Gamma, \\
\gamma_{N}^{\pm} \mathbf{u}^\pm(\mathbf{x}) &= \frac{1}{ik} \gamma_{D}^{\pm} \left( \nabla \times \mathbf{u}^\pm (\mathbf{x}) \right), \quad &\mathbf{x} \in \Gamma, 
\end{alignat}
along with their jumps and averages 
\begin{align}
\left[ \gamma_{\cdot} \right] := \gamma_{\cdot} ^+ - \gamma_{\cdot} ^-, \qquad
\left\{ \gamma_{\cdot} \right\} := \frac{1}{2} \left(\gamma_{\cdot} ^+ + \gamma_{\cdot} ^-\right).
\end{align}
For the potentials one has the trace and jump relations 
\begin{align}
\gamma ^\pm _{N} \mathcal{E} = \gamma ^\pm _{D} \mathcal{H}, \quad & \quad
\gamma ^\pm _{N} \mathcal{H} = -\gamma ^\pm _{D} \mathcal{E},\\
\left[\gamma _{D}\right] \mathcal{E}_\ell = \left[\gamma _{N}\right] \mathcal{H} = 0, \quad & \quad 
\left[\gamma _{N}\right] \mathcal{E} = \left[\gamma _{D}\right] \mathcal{H} = - \mathcal{I},
\end{align}
where $\mathcal{I}$ is the identity operator. 
We then define electric and magnetic boundary integral operators on $\Gamma$ by
\begin{alignat}{3}
&\mathcal{S} &:=& \{ \gamma_D \}\mathcal{E} &=& -\{ \gamma_{N} \} \mathcal{H},\\
&\mathcal{C} &:=& \{ \gamma _{D} \} \mathcal{H} &=& \{ \gamma _{N} \} \mathcal{E}, 
\end{alignat}
which by the jump relations also satisfy
\begin{alignat}{2}
\label{BIOrelns1}
\mathcal{S} &= \gamma_{D}^\pm \mathcal{E} = - \gamma_{N}^\pm \mathcal{H}, \\
\label{BIOrelns2}\mathcal{C} &= \gamma_{N}^\pm \mathcal{E} \pm \frac{1}{2} \mathcal{I} = \gamma_{D}^\pm \mathcal{H} \pm \frac{1}{2} \mathcal{I}.
\end{alignat}

By the Stratton-Chu formulae \cite{kirsch2015mathematical}, the interior and exterior fields $\mathbf{E}^i_m$, $m=1,\ldots,M$, and $\mathbf{E}^s$ in our multi-particle scattering problem can be represented as
\begin{alignat}{3}
\label{eqn:StrattonChu_int} \mathcal{H}^i_m (\gamma_{D,m}^- \mathbf{E}^i_m) &+ \mathcal{E}^i_m (\gamma_{N,m}^- \mathbf{E}^i_m) = 
\begin{cases}
\mathbf{E}^i_m(\mathbf{x}), & \mathbf{x} \in \Omega^i_m, \\
\mathbf{0}, & \mathbf{x} \not \in \overline{\Omega^i_m},
\end{cases}\\
\label{eqn:StrattonChu_ext}-\sum_m^M \mathcal{H}^e_m (\gamma_{D,m}^+ \mathbf{E}^s) &- \sum_m^M  \mathcal{E}^e_m (\gamma_{N,m}^+ \mathbf{E}^s) = 
\begin{cases}
\mathbf{E}^s(\mathbf{x}), & \mathbf{x} \in \Omega_e, \\
\mathbf{0}, & \mathbf{x} \not \in \overline{\Omega_e},
\end{cases}
\end{alignat}
where $(\mathcal{E}^i_m,\mathcal{H}^i_m,\gamma_{D,m}^-,\gamma_{N,m}^-)$ are $(\mathcal{E},\mathcal{H},\gamma_{D}^-,\gamma_{N}^-)$ for $\Gamma=\Gamma_m$ and $k=k_m$, and $(\mathcal{E}^e_m,\mathcal{H}^e_m,\gamma_{D,m}^+,\gamma_{N,m}^+)$ are $(\mathcal{E},\mathcal{H},\gamma_{D}^+,\gamma_{N}^+)$ for $\Gamma=\Gamma_m$ and $k=k_e$, for $m=1,\ldots, M$.

Taking appropriate interior and exterior Dirichlet and Neumann traces of (\ref{eqn:StrattonChu_int})-(\ref{eqn:StrattonChu_ext}), and recalling \eqref{BIOrelns1}-\eqref{BIOrelns2}, reveals that the boundary traces $(\gamma^-_{D,m}\mathbf{E}^i_m,\gamma^-_{N,m}\mathbf{E}^i_m, \gamma^+_{D,m}\mathbf{E}^s, \gamma^+_{N,m}\mathbf{E}^s)$ satisfy 
\begin{align}
\label{eqn:multiple_int}&\left( \frac{1}{2}\bm{\mathcal{I}}_m - \bm{\mathcal{{A}}}^i_m \right){\mathbf{u}}^i_m = 0, \\
\label{eqn:multiple_ext}&\left( \frac{1}{2} \bm{\mathcal{I}}_m + \bm{\mathcal{{A}}}^e_m \right) {\mathbf{u}}^s_m + \sum _{\ell \neq m}^M \bm{\mathcal{{A}}}_{m\ell} {\mathbf{u}}^s_\ell= 0 ,
\end{align}
where $\bm{\mathcal{I}}_m$ is the identity operator associated with $\Gamma_m$ and
\begin{gather}
\label{AimAem}
\bm{\mathcal{{A}}}^i_m = \begin{bmatrix}
\mathcal{C}^i_m & \frac{\mu_m}{k_m} \mathcal{S}^i_m \\[6pt]
-\frac{k_m}{\mu_m} \mathcal{S}^i_m & \mathcal{C}^i_m
\end{bmatrix}, \;
\bm{\mathcal{{A}}}^e_m = \begin{bmatrix}
\mathcal{C}^e_m & \frac{\mu_e}{k_e} \mathcal{S}^e_m \\[6pt]
-\frac{k_e}{\mu_e} \mathcal{S}^e_m & \mathcal{C}^e_m
\end{bmatrix}, \\
\label{Aml}
\bm{\mathcal{{A}}}_{m\ell} = \begin{bmatrix}
    \mathcal{C}^e_{m\ell} & \frac{\mu_e}{k_e} \mathcal{S}^e_{m\ell} \\[6pt]
    -\frac{k_e}{\mu_e} \mathcal{S}^e_{m\ell} & \mathcal{C}^e_{m\ell}
\end{bmatrix}, \\
{\mathbf{u}}^i_m = \begin{bmatrix}
    \gamma_{D,m}^{-} \mathbf{E}^i_m \\[6pt]
    \frac{k_m}{\mu_m} \gamma_{N,m}^- \mathbf{E}^i_m
\end{bmatrix}, \quad 
{\mathbf{u}}^s_m = \begin{bmatrix}
    \gamma_{D,m}^+ \mathbf{E}^s \\[6pt]
    \frac{k_e}{\mu_e} \gamma_{N,m}^+ \mathbf{E}^s
\end{bmatrix}.
\end{gather}
Here $(\mathcal{C}^i_m,\mathcal{S}^i_m)$ are $(\mathcal{C},\mathcal{S})$ for $\Gamma=\Gamma_m$ and $k=k_m$, and $(\mathcal{C}^e_m,\mathcal{S}^e_m)$ are $(\mathcal{C},\mathcal{S})$ for $\Gamma=\Gamma_m$ and $k=k_e$. The operators $(\mathcal{C}^e_{m\ell},\mathcal{S}^e_{m\ell})$ map from $\Gamma_\ell$ to $\Gamma_m$ and are defined for a boundary vector field $\mathbf{v}_\ell$ on $\Gamma_l$ by
\begin{align}
\label{}
\mathcal{C}^e_{m\ell}\mathbf{v}_\ell:= (\mathcal{C}^e \mathbf{\tilde{v}}_l)|_{\Gamma_m}
\quad
\mathcal{S}^e_{m\ell}\mathbf{v}_\ell:= (\mathcal{S}^e \mathbf{\tilde{v}}_l)|_{\Gamma_m}, 
\end{align}
where $\mathbf{\tilde{v}}_l$ denotes the vector field on $\cup_{j=1}^M \Gamma_j$ which equals $\mathbf{v}_\ell$ on $\Gamma_\ell$ and zero on $\Gamma_j$, for $j\neq \ell$; $|_{\Gamma_m}$ is restriction to $\Gamma_m$; and $(\mathcal{C}^e,\mathcal{S}^e)$ are $(\mathcal{C},\mathcal{S})$ for $\Gamma=\cup_{j=1}^M \Gamma_j$ and $k=k_e$.

The matrices $\left( \frac{1}{2}\bm{\mathcal{I}}_m + \bm{\mathcal{{A}}}^i_m \right) $ and $\left( \frac{1}{2}\bm{\mathcal{I}}_m - \bm{\mathcal{{A}}}^e_m \right)$
are scaled versions of the interior and exterior electromagnetic Calder\'on projectors on $\Gamma_m$, and satisfy the relations \cite{buffa2003galerkin}
\begin{align}
\label{eqn:calderon_int}\left( \frac{1}{2}\bm{\mathcal{I}}_m + \bm{\mathcal{{A}}}^i_m \right)^2&=\left( \frac{1}{2}\bm{\mathcal{I}}_m + \bm{\mathcal{{A}}}^i_m \right),\\
\label{eqn:calderon_ext}\left( \frac{1}{2}\bm{\mathcal{I}}_m - \bm{\mathcal{{A}}}^e_m \right)^2&=\left( \frac{1}{2}\bm{\mathcal{I}}_m - \bm{\mathcal{{A}}}^e_m \right).
\end{align}
These relations are central to the idea of Calder\'on preconditioning, as we shall explain shortly.
\section{The PMCHWT boundary integral formulation}\label{sct:PMCHWT}

Equation \eqref{eqn:multiple_int} is a system of boundary integral equations satisfied by the interior traces $(\gamma^-_{D,m}\mathbf{E}^i_m,\gamma^-_{N,m}\mathbf{E}^i_m)$ on $\Gamma_m$, and equation \eqref{eqn:multiple_ext} is a system of boundary integral equations satisfied by the exterior traces $(\gamma^+_{D,m}\mathbf{E}^s, \gamma^+_{N,m}\mathbf{E}^s)$. It is important to remark that they hold for any solutions of the Maxwell equations  \eqref{single_equation_i}-\eqref{single_equation_e}. To obtain the solution of our particular dieletric scattering problem we need to combine equations \eqref{eqn:multiple_int}-\eqref{eqn:multiple_ext} with the transmission conditions (\ref{bc1})-(\ref{bc2}), which we can rewrite as
\begin{alignat}{3}
\label{eqn:multiple_bc}{\mathbf{u}}^i_m = {\mathbf{u}}^s_m + {\mathbf{u}}^{inc}_m,\quad m=1,\ldots, M,
\end{alignat}
with
\begin{alignat}{3}
{\mathbf{u}}^{inc}_m &=& \begin{bmatrix}
    \gamma_{D,m}^+ \mathbf{E}^{inc} \\[6pt]
    \frac{k_e}{\mu_e} \gamma_{N,m}^+ \mathbf{E}^{inc}
\end{bmatrix}.
\end{alignat}
Equations (\ref{eqn:multiple_int})-(\ref{eqn:multiple_ext}) and \eqref{eqn:multiple_bc} can be combined in numerous different ways, leading to a range of different boundary integral equation formulations \cite{muller2013foundations, niino2012calderon}. Here we focus on the well-studied PMCHWT formulation \cite{poggio1970integral, wu1977scattering, mautz1977electromagnetic,harrington1989boundary}, which is obtained by subtracting (\ref{eqn:multiple_int}) from (\ref{eqn:multiple_ext}), then eliminating ${\mathbf{u}}^i_m$ using \eqref{eqn:multiple_bc}, to obtain, for each $m=1,\ldots,M$, the system
\begin{alignat}{1}
\left(\bm{\mathcal{{A}}}^i_m + \bm{\mathcal{{A}}}^e_m\right) {\mathbf{u}}^s_m + \sum _{\ell \neq m}^j \bm{\mathcal{{A}}}_{m\ell} {\mathbf{u}}^s_\ell = \left( \frac{1}{2}\bm{\mathcal{I}} - \bm{\mathcal{{A}}}^i_m \right) {\mathbf{u}}^{inc}_m.
\end{alignat}
We can combine these $M$ systems into a  block system
\begin{alignat}{3}
\label{eqn:multiple_scatterers}\bm{\mathcal{A}}\mathbf{ u}^s = \left(\frac{1}{2}\bm{\mathcal{ {I}}} - \bm{\mathcal{{A}}}^i \right) \mathbf{u}^{inc},
\end{alignat}
where 
\begin{alignat}{3}
\label{multiple_operator}&\bm{\mathcal{A}} &=& 
\begin{tikzpicture}[baseline={([yshift=-.5ex]current bounding box.center)}, ampersand replacement=\&]
\matrix (m) [matrix of math nodes,nodes in empty cells,right delimiter={]},left delimiter={[} ]{
\bm{\mathcal{{A}}}^e_1+\bm{\mathcal{{A}}}^i_1  \& \bm{\mathcal{{A}}}_{12}   \& \cdots  \& \bm{\mathcal{{A}}}_{1M}  \\
\bm{\mathcal{{A}}}_{21}    \& \& \& \vdots \\
\vdots  \&   \& \& \bm{\mathcal{{A}}}_{(M-1)M}   \\
\bm{\mathcal{{A}}}_{M1}  \& \cdots  \& \bm{\mathcal{{A}}}_{M(M-1)}   \& \bm{\mathcal{{A}}}^e_M+\bm{\mathcal{{A}}}^i_M\\
} ;
\draw[loosely dotted,thick] (m-1-1)-- (m-4-4);
\draw[loosely dotted,thick] (m-1-2)-- (m-3-4);
\draw[loosely dotted,thick] (m-2-1)-- (m-4-3);
\end{tikzpicture},
\end{alignat}
\begin{alignat}{3}
&\bm{\mathcal{A}}^i &=& 
\begin{tikzpicture}[baseline={([yshift=-.5ex]current bounding box.center)}, ampersand replacement=\&]
\matrix (m) [matrix of math nodes,nodes in empty cells,right delimiter={]},left delimiter={[} ]{
\bm{\mathcal{{A}}}^i_1  \& 0   \& \cdots  \& 0  \\
0    \& \& \& \vdots \\
\vdots  \&   \& \& 0   \\
0 \& \cdots  \& 0 \& \bm{\mathcal{{A}}}^i_M\\
} ;
\draw[loosely dotted,thick] (m-1-1)-- (m-4-4);
\draw[loosely dotted,thick] (m-1-2)-- (m-3-4);
\draw[loosely dotted,thick] (m-2-1)-- (m-4-3);
\end{tikzpicture},
\;
\bm{\mathcal{I}} = 
\begin{tikzpicture}[baseline={([yshift=-.5ex]current bounding box.center)}, ampersand replacement=\&]
\matrix (m) [matrix of math nodes,nodes in empty cells,right delimiter={]},left delimiter={[} ]{
\bm{\mathcal{{I}}}_1  \& 0   \& \cdots  \& 0  \\
0    \& \& \& \vdots \\
\vdots  \&   \& \& 0   \\
0 \& \cdots  \& 0 \& \bm{\mathcal{{I}}}_M\\
} ;
\draw[loosely dotted,thick] (m-1-1)-- (m-4-4);
\draw[loosely dotted,thick] (m-1-2)-- (m-3-4);
\draw[loosely dotted,thick] (m-2-1)-- (m-4-3);
\end{tikzpicture}, \\
%
&\mathbf{u}^s &=& 
\begin{tikzpicture}[baseline={([yshift=-.5ex]current bounding box.center)}, ampersand replacement=\&]
\matrix (m) [matrix of math nodes,nodes in empty cells,right delimiter={]},left delimiter={[} ]{
{\mathbf{u}}^s_1    \\
{\mathbf{u}}^s_2  \\
\vdots   \\
{\mathbf{u}}^s_M\\
};
\end{tikzpicture},  \quad
\mathbf{u}^{inc} = 
\begin{tikzpicture}[baseline={([yshift=-.5ex]current bounding box.center)}, ampersand replacement=\&]
\matrix (m) [matrix of math nodes,nodes in empty cells,right delimiter={]},left delimiter={[} ]{
{\mathbf{u}}^{inc}_1    \\
{\mathbf{u}}^{inc}_2  \\
\vdots   \\
{\mathbf{u}}^{inc}_M\\
};  
\end{tikzpicture}.
\end{alignat}
Equation \eqref{eqn:multiple_scatterers} is the 
PMCHWT formulation, expressed in multi-particle notation. 

\section{Galerkin method and preconditioning}\label{sct:preconditioning}
Equation \eqref{eqn:multiple_scatterers} is known to be well-posed, when we view $\bm{\mathcal{A}}$ as a mapping $\bm{\mathcal{A}}:X\to X$ on the function space \mbox{$X=\oplus_{m=1}^M \mathbf{H}^{-\frac{1}{2}}_\times (\textnormal{div}_{\Gamma_m}, \Gamma_m)^2$}, where $\mathbf{H}^{-\frac{1}{2}}_\times (\textnormal{div}_{\Gamma_m}, \Gamma_m)$ denotes the space of tangential vector fields on $\Gamma_m$ of Sobolev regularity $-1/2$ whose surface divergences  also have Sobolev regularity $-1/2$ (see \cite{buffa2003galerkin} for details). 
Upon Galerkin discretization (described in more detail below) 
it yields a system of linear equations, the solution of which can be inserted into the representation formulae \eqref{eqn:StrattonChu_int}-\eqref{eqn:StrattonChu_ext} to produce a solution of the original scattering problem. 

Unfortunately, the resulting linear system is ill-conditioned, leading to slow convergence of iterative solvers such as GMRES. 
The origin of the ill-conditioning can be understood from the properties of the underlying continuous operators. It is enough to consider the single-particle scattering case ($M=1$), in which case the operator $\bm{\mathcal{{A}}}$ is
\begin{alignat}{3}
\bm{\mathcal{{A}}}=\bm{\mathcal{{A}}}^e_1 + \bm{\mathcal{{A}}}^i_1 = \begin{bmatrix}
\mathcal{C}^e_1 + \mathcal{C}^i_1 & \frac{\mu_e}{k_e}\mathcal{S}^e_1 + \frac{\mu_1}{k_1}\mathcal{S}^i_1 \\[6pt]
-\frac{k_e}{\mu_e} \mathcal{S}^e_1 -\frac{k_1}{\mu_1} \mathcal{S}^i_1 & \mathcal{C}^e_1 + \mathcal{C}^i_1
\end{bmatrix}.
\end{alignat}
Ill-conditioning at the discrete level should be expected because the operators $\mathcal{C}^i_1$ and $\mathcal{C}^e_1$ are compact \cite{nedelec2001acoustic}, with eigenvalues accumulating at zero, while the operators $\mathcal{S}^i_1$ and $\mathcal{S}^e_1$ are both the sum of a compact operator, with eigenvalues accumulating at zero, and a hypersingular operator, with eigenvalues accumulating at infinity. For a more detailed discussion we refer the reader to \cite{cools2011calderon}.

For the single-particle scattering problem it was shown in \cite{yan2010comparative, cools2011calderon, niino2012calderon} that the conditioning of the linear system can be dramatically improved by ``Calder\'on preconditioning'', a form of operator preconditioning (i.e., applied at the continuous level, before discretization) that exploits the projection properties \eqref{eqn:calderon_int}-\eqref{eqn:calderon_ext} of the Calder\'on projectors. 
Explicitly, it involves applying the operator  $\bm{\mathcal{{A}}}$ to both sides of (\ref{eqn:multiple_scatterers}) to give the equivalent equation
\begin{alignat}{3}
\label{squared_system}
\bm{\mathcal{{A}}}^2 {\mathbf{u}}^s =\bm{\mathcal{{A}}}  \left( \frac{1}{2}\bm{\mathcal{I}} - \bm{\mathcal{{A}}}^i \right) {\mathbf{u}}^{inc}.
\end{alignat}
Squaring the operator in this way regularizes the system by shifting the accumulation points of the spectrum away from the origin and taming the hypersingular component. For details we refer to \cite{yan2010comparative, cools2011calderon, niino2012calderon}, but the regularization relies on the following relations (contained in \eqref{eqn:calderon_int}-\eqref{eqn:calderon_ext}):
\begin{align}
\label{calderonproperty1}(\mathcal{S}^i_1)^2  = -\frac{1}{4} \mathcal{I}_1 + (\mathcal{C}^i_1)^2, \quad & \quad (\mathcal{S}^e_1)^2  = -\frac{1}{4} \mathcal{I}_1 + (\mathcal{C}^e_1)^2, \\
\label{calderonproperty2}\mathcal{C}^i_1 \mathcal{S}^i_1 + \mathcal{S}^i_1 \mathcal{C}^i_1 = 0, \quad & \quad  \mathcal{C}^e_1 \mathcal{S}^e_1 + \mathcal{S}^e_1 \mathcal{C}^e_1 = 0.
\end{align}
In particular, (\ref{calderonproperty1}) implies that $(\mathcal{S}^i_1)^2$ and $(\mathcal{S}^e_1)^2$ are second-kind integral operators (i.e., of the form ``constant times identity plus compact operator'') with eigenvalues accumulating at $-1/4$. This is referred to as the ``self-regularising property'' of the operators $\mathcal{S}^i_1$, $\mathcal{S}^e_1$ \cite{yan2010comparative}.

For the multi-particle scattering case ($M>1$) one can adopt the same preconditioning strategy and solve the squared system \eqref{squared_system} instead of (\ref{eqn:multiple_scatterers}). However, since it is the diagonal blocks in \eqref{multiple_operator} that cause the ill-conditioning (the off-diagonal blocks are compact since they map between different boundary components), one might imagine that it is sufficient to precondition block-diagonally by
\begin{alignat}{3}
\label{eqn:diagonal_preconditioner}\bm{\mathcal{D}} \bm{\mathcal{{A}}} {\mathbf{u}}^s = \bm{\mathcal{D}} \left( \frac{1}{2}\bm{\mathcal{I}} - \bm{\mathcal{{A}}}^i \right) {\mathbf{u}}^{inc},
\end{alignat}
where 
\begin{align}
\label{}
\bm{\mathcal{D}} = 
\begin{tikzpicture}[baseline={([yshift=-.5ex]current bounding box.center)}, ampersand replacement=\&]
\matrix (m) [matrix of math nodes,nodes in empty cells,right delimiter={]},left delimiter={[} ]{
\bm{\mathcal{{A}}}^e_1+\bm{\mathcal{{A}}}^i_1  \& 0   \& \cdots  \& 0  \\
0    \& \& \& \vdots \\
\vdots  \&   \& \& 0   \\
0 \& \cdots  \& 0 \& \bm{\mathcal{{A}}}^e_M+\bm{\mathcal{{A}}}^i_M\\
} ;
\draw[loosely dotted,thick] (m-1-1)-- (m-4-4);
\draw[loosely dotted,thick] (m-1-2)-- (m-3-4);
\draw[loosely dotted,thick] (m-2-1)-- (m-4-3);
\end{tikzpicture}.
\end{align}
As we will show experimentally, this block-diagonal preconditioner achieves a similar improvement in conditioning as for \eqref{squared_system} but with a reduced computational cost. 

Implementing Calder\'on preconditioners requires the discretization of products of boundary integral operators. This does not simply involve taking the product of the respective Galerkin matrices. Instead one needs to introduce appropriate mass matrices, which arise in the strong form of the Galerkin operators. 
For further details about weak and strong discrete forms we refer to \cite{betcke2017product}; here we just outline the details relevant to the problem at hand.

We start by noting that since
$\mathbf{H}^{-\frac{1}{2}}_\times (\textnormal{div}_{\Gamma_m}, \Gamma_m)$ is self-dual with respect to the twisted $L^2$ dual pairing \cite{buffa2003galerkin}
\[ \langle \mathbf{a},\mathbf{b} \rangle_{\Gamma_m} = \int_{\Gamma_m} \mathbf{a}\cdot (\mathbf{n}\times \mathbf{b})\,\mathrm{d}S, \]
the space $X$ is self-dual with respect to the pairing
\begin{align}
\label{dualdef}
&\left\langle \bigoplus_{m=1}^M\left(\begin{array}{c}\mathbf{c}_m \\ \mathbf{d}_m \end{array}\right),\bigoplus_{m=1}^M\left(\begin{array}{c}\mathbf{e}_m \\ \mathbf{f}_m \end{array}\right)\right\rangle\notag
\\
&\qquad := \sum_{m=1}^M \langle \mathbf{c}_m ,\mathbf{f}_m \rangle_{\Gamma_m} + \langle\mathbf{d}_m , \mathbf{e}_m \rangle_{\Gamma_m}.
\end{align}
To define a Galerkin method we choose discrete trial and test spaces $X_h,Y_h\subset X$ of some common dimension $J\in \mathbb{N}$, with bases $\{\bm{\psi}_j\}_{j=1}^J$ and $\{\bm{\phi}_j\}_{j=1}^J$. The Galerkin solution $\mathbf{u}^s_h = \sum_{j=1}^J x_j \bm{\psi}_j\in X_h$ satisfies the variational problem
\begin{align}
\label{}
\langle \bm{\mathcal{A}}\mathbf{u}^s_h,\mathbf{v}_h \rangle = \left< \left(\frac{1}{2}\bm{\mathcal{ {I}}} - \bm{\mathcal{{A}}}^i \right) \mathbf{u}^{inc},\mathbf{v}_h \right>, \quad \forall \mathbf{v}_h  \in Y_h,
\end{align}
which corresponds to the linear system
\begin{align}
\label{linsys}
\mathbf{A}\mathbf{x} = \mathbf{b},
\end{align}
where 
$\mathbf{x} = (x_1,\ldots,x_J)^T$, 
$\mathbf{b} = (b_1,\ldots,b_J)^T$, with 
$b_j=\langle \left(\frac{1}{2}\bm{\mathcal{ {I}}} - \bm{\mathcal{{A}}}^i \right) \mathbf{u}^{inc},\bm{\phi}_j \rangle$
and $\mathbf{A}$ is the Galerkin matrix with 
\begin{align}
\label{}
\mathbf{A}_{ij} = \langle \bm{\mathcal{A}}\bm{\psi}_j,\bm{\phi}_i \rangle, \quad i,j=1,\ldots,J.
\end{align}
This matrix corresponds to the discrete weak form of $\bm{\mathcal{A}}$, which 
maps $X_h$ to the dual space of $Y_h$.  
The discrete strong form of $\bm{\mathcal{A}}$, which maps $X_h$ to $X_h$, has matrix 
$\mathbf{M}^{-1}\mathbf{A}$, 
where $\mathbf{M}$ is the mass matrix with entries
\begin{align}
\label{}
\mathbf{M}_{ij} = \langle \bm{\phi}_i,\bm{\psi}_j \rangle, \quad i,j=1,\ldots,J.
\end{align}
The conversion of \eqref{linsys} to the system
\begin{align}
 \label{}
\mathbf{M}^{-1}\mathbf{A}\mathbf{x} = \mathbf{M}^{-1}\mathbf{b}
 \end{align} 
is sometimes known as ``mass-matrix preconditioning''. 
The weak and strong forms of the squared operator $\bm{\mathcal{A}}^2$ have matrices $\mathbf{A}\mathbf{M}^{-1}\mathbf{A}$ and $\mathbf{M}^{-1}\mathbf{A}\mathbf{M}^{-1}\mathbf{A}$ respectively, and those of the block-diagonally preconditioned operator $\bm{\mathcal{D}}\bm{\mathcal{A}}$ have matrices $\mathbf{D}\mathbf{M}^{-1}\mathbf{A}$ and $\mathbf{M}^{-1}\mathbf{D}\mathbf{M}^{-1}\mathbf{A}$ respectively, where $\mathbf{D}$ is the Galerkin matrix for $\bm{\mathcal{D}}$. Assuming that the basis functions are ``local'', and are indexed in a natural way so that there exists $1=J_1<J_2<\ldots <J_M<J_{M+1}=J+1$ with ${\rm supp}\,\bm{\phi}_j,{\rm supp}\,\bm{\psi}_j \subset \Gamma_m$ for $J_m\leq j<J_{m+1}$, the matrix $\mathbf{M}$ is block-diagonal and $\mathbf{D}$ is the block-diagonal part of $\mathbf{A}$.

For $\mathbf{M}$ to be invertible we need $X_h$ to be dual to $Y_h$ with respect to the pairing $\langle \cdot,\cdot\rangle$, i.e.\ for $\langle \cdot,\cdot\rangle$ to be inf-sup stable on $X_h\times Y_h$. This places a restriction on the possible choices of trial and test spaces $X_h$ and $Y_h$. In our numerical experiments we use
\begin{align}
\label{}
X_h=Y_h=\bigoplus_{m=1}^M \left(\begin{array}{c}{\rm RWG}_m \\ {\rm BC}_m \end{array}\right),
\end{align}
where ${\rm RWG}_m$ denotes the span of the Rao-Wilton-Glisson basis functions \cite{rwg1982} on the primal mesh of $\Gamma_m$ and $\rm BC$ denotes the span of the Buffa-Cristiansen basis functions \cite{buffa2007dual} on the barycentric dual mesh of $\Gamma_m$. 
Here we are using the fact that ${\rm BC}_m$ and ${\rm RWG}_m$ are dual with respect to the twisted duality pairing $\langle \cdot,\cdot\rangle_{\Gamma_m}$.
Using different discretizations for the two components, combined with the dual pairing \eqref{dualdef}, is attractive because it produces a symmetric formulation, i.e.\ the trial and test spaces coincide. 
For more details of this formulation and its implementation in Bempp,
see \cite{scroggs2017software}. In particular we note that in the Bempp implementation, the
BC basic functions are ordered and rotated consistently with their RWG
siblings.
But we note that other discretizations are possible \cite{cools2011calderon,niino2012calderon}
and are just as easily implemented with the Bempp software package. 

\section{Computational complexity}\label{sct:complexity}
\begin{table*}[!t]
\centering
\begin{tabular}{lll}
\toprule
Continuous operator  & Discrete Weak Form  & Discrete Strong Form        \\
\midrule
$\bm{\mathcal{A}}$ & $\mathbf{A}\mathbf{x} = \mathbf{b}$& $\mathbf{M}^{-1}\mathbf{A}\mathbf{x} = \mathbf{M}^{-1}\mathbf{b}$ \\[5pt]
$\bm{\mathcal{A}}^2$ & $\mathbf{A}\mathbf{M}^{-1}\mathbf{A}\mathbf{x} = \mathbf{A}\mathbf{M}^{-1}\mathbf{b}
$ & $\mathbf{M}^{-1}\mathbf{A}\mathbf{M}^{-1}\mathbf{A}\mathbf{x} = \mathbf{M}^{-1}\mathbf{A}\mathbf{M}^{-1}\mathbf{b}
$\\[5pt]
$\bm{\mathcal{D}}\bm{\mathcal{A}}$ & $\mathbf{D}\mathbf{M}^{-1}\mathbf{A}\mathbf{x} = \mathbf{D}\mathbf{M}^{-1}\mathbf{b}
$ & $\mathbf{M}^{-1}\mathbf{D}\mathbf{M}^{-1}\mathbf{A}\mathbf{x} = \mathbf{M}^{-1}\mathbf{D}\mathbf{M}^{-1}\mathbf{b}
$\\[5pt]
\bottomrule
\end{tabular}
\caption{Linear systems for the weak and strong discrete forms of the continuous operators.}
\label{table:discrete_operators}
\end{table*}

For multi-particle scattering problems we have introduced six different formulations involving the PMCHWT operator $\bm{\mathcal{A}}$ (in weak or strong form), the squared operator $\bm{\mathcal{A}}^2$ (in weak or strong form), and the block-diagonally preconditioned operator $\bm{\mathcal{D}}\bm{\mathcal{A}}$ (in weak or strong form). For reference we present in  Table \ref{table:discrete_operators} the linear systems that have to be solved in each of these six formulations.

We solve these linear systems using the iterative method GMRES. The performance of each preconditioning strategy will be assessed by measuring the computational cost required to solve the linear system to within a certain prescribed tolerance on the relative residual.
This computational cost depends on both the number of GMRES iterations required, and the cost of each iteration.

To measure the cost of each strategy we will count the total number of matrix-vector products (termed ``matvecs'' henceforth) incurred. By a single matvec we mean a single application of one discretized boundary integral operator $\mathcal{C}^i_m$, $\mathcal{C}^e_m$, $\mathcal{S}^i_m$, $\mathcal{S}^e_m$ etc. 
Applications of the inverse mass matrix $\mathbf{M}^{-1}$ are not included in the matvec count, since their cost is negligible compared to those of the other operators. 
In more detail, the matrix $\mathbf{M}$ is block-diagonal with $M$ sparse blocks. Computing its LU decomposition requires $M$ independent sparse LU decompositions (one for each scatterer), which can be effectively parallelized. For two-dimensional problems the computation of the LU decomposition of a mass matrix is very efficient using suitable reordering strategies of the elements. In our case the mass-matrix is formed over a two dimensional manifold in three dimensional space. This makes the LU decomposition slightly more expensive, but still reasonably cheap, given that it only needs to be performed once during a precomputation for each domain.

Recalling \eqref{multiple_operator} and \eqref{AimAem}-\eqref{Aml}, we note that a single application of the matrix $\mathbf{A}$ requires $4M(M+1)$ matvecs ($4$ for each off-diagonal block and $8$ for each diagonal block), and that a single application of the matrix $\mathbf{D}$ requires $8M$ matvecs. Hence the overall cost of the three formulations (taking into account the initial pre-multiplication of the right-hand-side) is
\begin{align}
\label{matvecA}
\bm{\mathcal{A}}
&: 4M(M+1)(G + \floor*{ G/\rho}) \text{ matvecs}, \\
\label{matvecA2}
\bm{\mathcal{A}}^2     &: 8M(M+1)(G + \floor*{ G/\rho}) + 4M(M+1)\text{ matvecs}, \\
\label{matvecDA}
\bm{\mathcal{D}}\bm{\mathcal{A}}
&: 4M(M +3) (G + \floor*{ G/\rho}) + 8M\text{ matvecs},
\end{align}
where $G$ is the number of GMRES iterations required to achieve the specified tolerance, $\rho$ is the number of iterations per GMRES cycle passed as the \texttt{restart} argument in GMRES, and $\floor*{\cdot}$ is the ``floor'' function.  Note that while the matvec count per GMRES iteration is the same for the weak and strong forms of each formulation (because we are excluding mass matrix solves, as explained above), the value of $G$ (and hence the overall matvec count) will in general be very different for the weak and strong forms. Indeed, the advantage of working with strong forms is one of the key messages of the paper.

Of course the overall computational cost is not governed simply by the total matvec count, since it also depends on the discretization resolution (mesh width) and the method used to assemble the operators. We now make some brief remarks on these matters. 
Regarding mesh size, suppose that for each $m=1,\ldots,M$, the scatterer $\Gamma_m$ is discretized on a triangular mesh with $N_m$ elements, and set $N=\max_m N_m$. In order to capture the wave solution one needs to use a fixed number (typically around $10$) of elements per exterior wavelength $\lambda_e := 2\pi / k_e$, so that $N\sim k_e^2$. 
Regarding assembly, for small problems a dense matrix discretization of the operators is possible, which results in a computational complexity for the assembly and matvec of $\mathcal{O}(N^2)$, or equivalently $\mathcal{O}(k_e^4)$. For practical applications this is too expensive. Alternatives are hierarchical matrix ($\mathcal{H}$-Matrix) \cite{hackbusch2015hierarchical} or fast multipole methods (FMM) \cite{gumerov2005fast}. $\mathcal{H}$-Matrices are most effective for problems with only a moderate number of wavelengths. Their complexity is $\mathcal{O}(rN\log N)$ for assembly and matvec, where $r$ is a measure of the local approximation rank required to achieve the prescribed accuracy. For non-oscillatory problems $r$ is effectively constant. For high-frequency problems we asymptotically have $r\sim k_e^2\sim N$ \cite{engquist2018approximate}. However, this is a worst-case estimate and in practical applications one often observes a complexity of $\mathcal{O}(N^{\alpha}\log N)$ for some $1<\alpha<2$, with $\alpha$ typically being close to $1$ even for highly oscillatory problems \cite{betcke2017computationally}. 
If instead of $\mathcal{H}$-Matrices a high-frequency FMM implementation is chosen then the complexity is the same but with $\alpha=1$. The price to pay is that a stable high-frequency FMM is significantly more challenging to implement than a standard $\mathcal{H}$-Matrix method. Both methods also differ in their realistic timing behaviour for assembly and matvecs. As a rule of thumb, while $\mathcal{H}$-Matrices have a longer assembly time and faster matvecs, classical high-frequency FMM has a shorter assembly time but slower matvecs. Hence, optimising the number of matvecs is especially important for FMM and high-frequency applications. It is also vital when one needs to perform multiple solves for different right-hand-sides, e.g.\ when computing averages over particle orientation \cite{groth2015boundary}.

The numerical results in this paper were computed using the Bempp boundary element library (available at \url{www.bempp.com}). It has a built-in, custom-developed $\mathcal{H}$-Matrix implementation, which offers thread-based parallelization on single nodes and MPI based parallelization on distributed nodes. Sparse Matrix LU factorizations (for the mass matrices) are computed through the SuperLU solver interfaced in the sparse matrix module of Python's scipy library (\url{www.scipy.org}). The SuperLU version used in Scipy does not offer parallelization, and moreover, Bempp by default cannot parallelize the block-diagonal mass matrix computation across the scatterers in a multi-scattering configuration. 
While the current Bempp version offers sufficient performance for small to medium sized scattering configurations, a scalable high-frequency FMM implementation is currently in development.

\section{Benchmarks}\label{sct:benchmarks}

In this section we compare the performance of the six discrete formulations in Table \ref{table:discrete_operators} on a range of benchmark problems. 
In all our experiments the incident wave is a plane wave, and the maximum BEM mesh element size is $2 \pi /(10 k_e)$, so that there are at least 10 elements per exterior wavelength in each coordinate direction. As was investigated in detail in \cite{groth2015boundary}, for the type of dielectric scattering problems considered here, this leads to a typical discretization error of approximately 1\% ($10^{-2}$ relative error). (In particular, we note that, since all our scatterers are polyhedral, the scatterer geometry is captured exactly by the BEM mesh, the only exception being the spheres in Figure 5.)
In all our experiments we terminate the GMRES solver once the relative residual falls below $10^{-5}$, with the \texttt{restart} option having the default value of $\rho = 20$ iterations.

We present results for two different refractive indices: one with weak absorption $n_1=1.311 + 2.289 \times 10^{-9}\mathrm{i}$  and one with high absorption $n_2 = 1.0833 + 0.204 \mathrm{i}$. These correspond to the measured refractive index of ice at the wavelengths $\lambda_1 = 0.55\mu m$ and $\lambda_2 = 10.87\mu m$ respectively \cite{warren2008optical}. Using these two representative refractive indices we present numerical results for scattering by different configurations of particles of fixed nondimensional size at a range of wavenumbers, simulating scattering by fixed geometrical configurations at a range of different size parameters. Specifically, as our definition of size parameter we use the quantity $k_er$, where $r$ is the radius of the smallest sphere enclosing the entire scatterer configuration.

\begin{table*}[!t]
\centering
\begin{tabular}{lrrrrrrr}
\toprule
& \multicolumn{3}{c}{$n=1.311 + 2.289 \times 10^{-9}\mathrm{i}$} & \phantom{abc} & \multicolumn{3}{c}{$n=1.0833 + 0.204 \mathrm{i}$} \\
\cmidrule{2-4} \cmidrule{6-8} 
& $k_e=4$  & $k_e=6$  & $k_e=10$  & & $k_e=4$  & $k_e=6$  & $k_e=10$      \\
\midrule
Discrete operator   &   &   \\
$\mathbf{A}$ & 599 (5024) &   548 (4600) &  270 (2264) & & 332 (2784) &   316 (2648) &  366 (3072)   \\[4pt]
$\mathbf{M}^{-1} \mathbf{A}$ & 11 (88) & 13 (104) & 18 (144)  && 9 (72) & 9 (72) & 10 (80)\\[4pt]
$\mathbf{A}\mathbf{M}^{-1}\mathbf{A}$ &34 (568) &38 (632)   & 58 (968) &&30 (504) &31 (520)  &34 (568)    \\[4pt]
$\mathbf{M}^{-1}\mathbf{A}\mathbf{M}^{-1}\mathbf{A}$ &6 (104) &  7 (120) &   10 (168) &&5 (88) &  5 (88) & 5 (88)   \\
\bottomrule
\end{tabular}
\caption{Number of GMRES iterations and total matvec count (in brackets) for the different discrete formulations for scattering by a single unit cube ($M=1$). The mesh size is $h=2\pi/(10 k_e)$, and the incident wave is $\mathbf{E}^{inc}(\mathbf{x})=\mathbf{p} \mathrm{e}^{\mathrm{i} k_e \mathbf{d} \cdot \mathbf{x}}$, with $\mathbf{d} = (1,0,0)^T$ and $\mathbf{p} = (0,0,1)^T$. The magnetic permeabilities are $\mu _1 = \mu _e = 1$. The size parameter $k_e r = \sqrt{3}/2 k_e$, with $k_e$  as defined in the three different cases.}
\label{table:iterations}
\end{table*}

\begin{figure*}[!t]
	\centering
    \includegraphics[width = 0.44\textwidth]{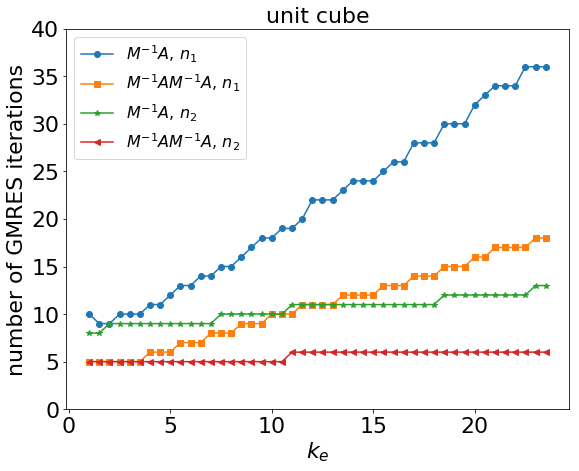}
    \includegraphics[width = 0.45 \textwidth]{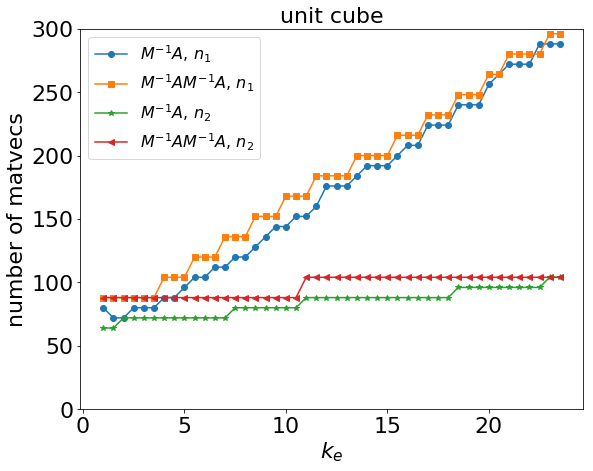}
\caption{Performance of the discrete strong form operators in the case of scattering by a single unit cube ($M=1$) as a function of exterior wavenumber $k_e$. Results are shown for  low absorption (refractive index $n_1=1.311 + 2.289 \times 10^{-9}\mathrm{i}$) and high absorption (refractive index $n_2=1.0833 + 0.204 \mathrm{i}$).
Other parameters are as in Table \ref{table:iterations}. }
\label{fig:iterations}
\end{figure*}

We begin with the case of single-particle scattering ($M=1$). 
In Table \ref{table:iterations} we present GMRES iteration and matvec counts for the weak and strong forms of the unpreconditioned operator $\bm{\mathcal{A}}$ and the Calder\'on-preconditioned operator $\bm{\mathcal{A}}^2$, for scattering by a unit cube. (Note that when $M=1$ it holds that $\bm{\mathcal{D}}\bm{\mathcal{A}}=\bm{\mathcal{A}}^2$.)
The ill-conditioning of the unpreconditioned weak form $\mathbf{A}$ is clearly visible in the large number of GMRES iterations (and hence matvecs) in the first row of the Table \ref{table:iterations}. The results in the lower rows of the table show that the other three formulations all provide a significant improvement in performance compared to the unpreconditioned weak form. For this single-particle scattering problem the unpreconditioned strong form $\mathbf{M}^{-1}\mathbf{A}$ (i.e.\ mass-matrix preconditioning) performs the best in terms of overall matvec count. While the Calder\'on preconditioned strong form $\mathbf{M}^{-1}\mathbf{A}\mathbf{M}^{-1}\mathbf{A}$ requires the lowest number of GMRES iterations, the increased matvec count for each iteration means that overall it is more expensive than $\mathbf{M}^{-1}\mathbf{A}$. 
While the Calder\'on preconditioned weak form $\mathbf{A}\mathbf{M}^{-1}\mathbf{A}$ (which is the formulation studied in \cite{cools2011calderon}) provides a significant improvement over the unpreconditioned weak form $\mathbf{A}$, it cannot compete with the strong forms $\mathbf{M}^{-1}\mathbf{A}\mathbf{M}^{-1}\mathbf{A}$ and $\mathbf{M}^{-1}\mathbf{A}$.

In Figure \ref{fig:iterations} we compare the performance of the two strong forms $\mathbf{M}^{-1}\mathbf{A}\mathbf{M}^{-1}\mathbf{A}$ and $\mathbf{M}^{-1}\mathbf{A}$ as a function of the exterior wavenumber $k_e$. 
For both refractive indices the number of GMRES iterations for $\mathbf{M}^{-1}\mathbf{A}\mathbf{M}^{-1}\mathbf{A}$ is approximately half that for $\mathbf{M}^{-1}\mathbf{A}$ (with faster convergence for high absorption), but the matvec count per iteration is more than doubled (see \eqref{matvecA}-\eqref{matvecA2}), so that the overall matvec count is lower for $\mathbf{M}^{-1}\mathbf{A}$.

\begin{figure*}[!t]
    \centering
    \includegraphics[width = 0.31 \textwidth]{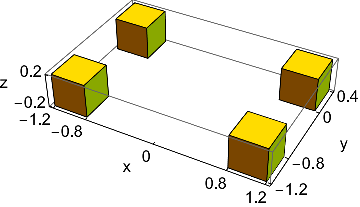}
    \hfill
    \includegraphics[width = 0.31 \textwidth]{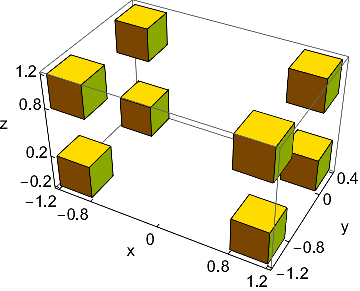}
    \hfill
    \includegraphics[width = 0.31 \textwidth]{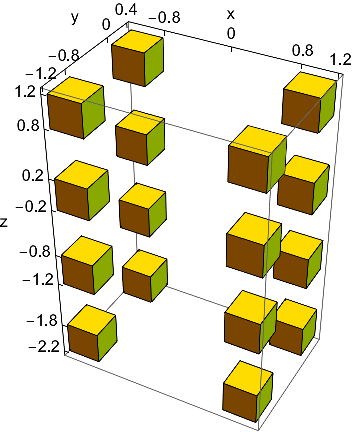}
    \caption{Arrangement of 4, 8 and 16 cubes of side length 0.4. The size parameters $k_e r$ are $1.4 k_e$, $1.6 k_e$ and $2.2k_e$ respectively, where  $k_e$ is the exterior wavenumber.}
    \label{fig:cubes}
\end{figure*}

We now consider scattering by multiple particles ($M>1$). In Table \ref{table:iterations_multiple} we present GMRES iteration and matvec counts for the weak and strong forms of the unpreconditioned operator $\bm{\mathcal{A}}$, the Calder\'on-preconditioned operator $\bm{\mathcal{A}}^2$, and the block-diagonally preconditioned operator $\bm{\mathcal{D}}\bm{\mathcal{A}}$, for scattering by an array of four cubes aligned as in Figure \ref{fig:cubes}. 
As in the single-particle scattering case, the strong forms perform significantly better than their corresponding weak forms. And again, mass matrix preconditioning $\mathbf{M}^{-1} \mathbf{A}$ performs better in terms of overall matvec count than Calder\'on strong form $\mathbf{M}^{-1} \mathbf{A}\mathbf{M}^{-1} \mathbf{A}$. But better than both of these methods in terms of overall matvec count is the strong form of the block-diagonally preconditioned formulation $\mathbf{M}^{-1} \mathbf{D}\mathbf{M}^{-1} \mathbf{A}$, which offers a similar GMRES iteration count to the full Calder\'on strong form $\mathbf{M}^{-1} \mathbf{A}\mathbf{M}^{-1} \mathbf{A}$, but at a significantly lower cost per iteration (see \eqref{matvecA}-\eqref{matvecDA}). 

\begin{table*}[!t]
\centering
\begin{tabular}{lrrrrrrrr}
\toprule
& \multicolumn{3}{c}{$n=1.311 + 2.289 \times 10^{-9}\mathrm{i}$} & \phantom{abc} & \multicolumn{3}{c}{$n=1.0833 + 0.204 \mathrm{i}$} \\
\cmidrule{2-4} \cmidrule{6-8} 
 & $k_e =7$  & $k_e =12$ & $k_e =22$ & & $k_e =7$  & $k_e =12$ & $k_e =22$    \\
\midrule
Discrete operator   &   &  & \\
$\mathbf{A}$  &326 (27360) &1193 (100160) &406 (34080) &  &190 (15920) &535 (44880) &382 (32080) \\
$\mathbf{M}^{-1} \mathbf{A}$  &12 (960) &15 (1200) &24 (2000) & &10 (800)  &11 (880) &12 (960) \\
$\mathbf{A}\mathbf{M}^{-1} \mathbf{A}$   &32 (5360)  &48 (8080) &122 (20560) & &30 (5040) &36 (6000) &42 (7120)  \\
$\mathbf{M}^{-1} \mathbf{A}\mathbf{M}^{-1} \mathbf{A}$  &6 (1040) &8 (1360) &12 (2000) & &5 (880) &6 (1040) &6 (1040) \\
$\mathbf{D}\mathbf{M}^{-1} \mathbf{A}$   &35 (4064) &48 (5632)  &69 (8096) & &36 (4176) &42 (4960) &44 (5184)\\
$\mathbf{M}^{-1} \mathbf{D}\mathbf{M}^{-1} \mathbf{A}$  &7 (816) &8 (928) &12 (1376) & &6 (704) &6 (704) &7 (816)  \\
\bottomrule
\end{tabular}
\caption{
Number of GMRES iterations and total matvec count (in brackets) for the different discrete formulations for scattering by an array of four cubes ($M=4$). The cubes have side length $0.4$ and are arranged as in Figure \ref{fig:cubes}. The mesh size is $h=2\pi/(10 k_e)$, and the incident wave is $\mathbf{E}^{inc}(\mathbf{x})=\mathbf{p} \mathrm{e}^{\mathrm{i} k_e \mathbf{d} \cdot \mathbf{x}}$, with $\mathbf{d} = (1/\sqrt{2},1/\sqrt{2},0)$ and $\mathbf{p} = (0,0,1)^T$. The magnetic permeabilities are $\mu _m = \mu _e = 1$, for $m=1, \ldots, 4$.}
\label{table:iterations_multiple}
\end{table*}

\begin{figure*}[t!]
\centering
    \begin{subfigure}[t]{\textwidth}
    \centering
        \includegraphics[width = 0.42 \textwidth]{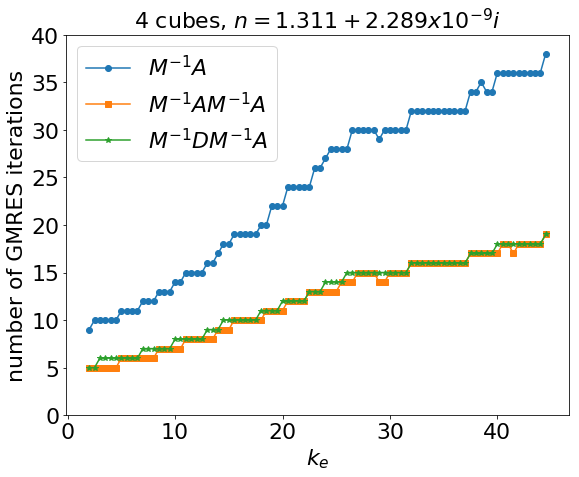}
        \includegraphics[width = 0.44 \textwidth]{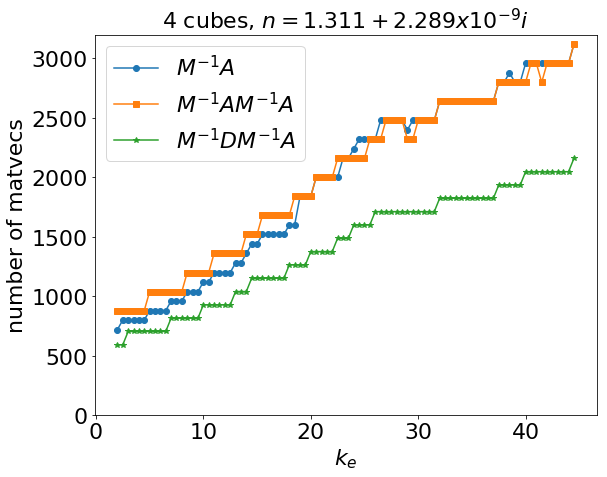}
        \caption{Low absorption}
    \end{subfigure}
    \vspace*{0.2cm}
    \begin{subfigure}[t]{\textwidth}
    \centering 
        \includegraphics[width = 0.42 \textwidth]{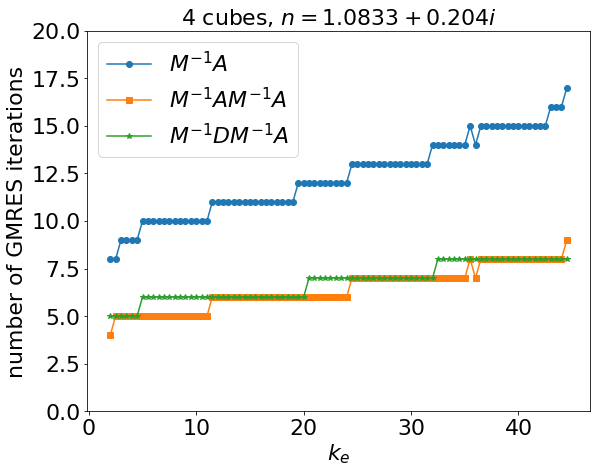}
        \includegraphics[width = 0.44 \textwidth]{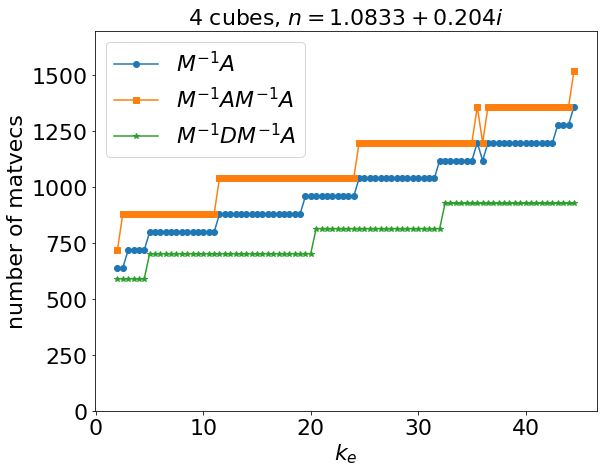}
        \caption{High absorption}
    \end{subfigure}
    \begin{subfigure}[t]{\textwidth}
    \centering
        \includegraphics[width = 0.42 \textwidth]{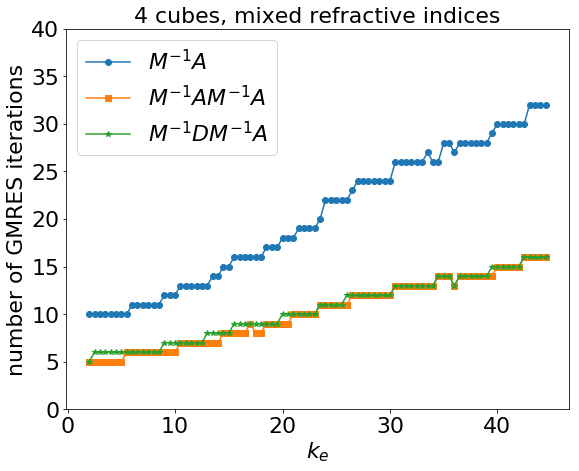}
        \includegraphics[width = 0.44 \textwidth]{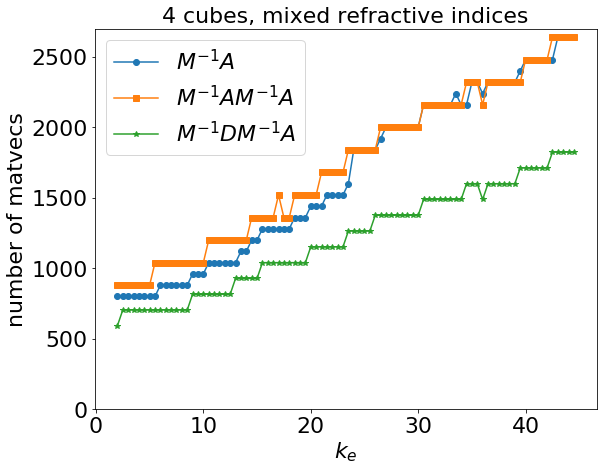}
        \caption{Mixed zero/low/high absorption}
    \end{subfigure}
\caption{
Performance of the discrete strong form operators in the case of scattering by an array of four cubes ($M=4$) as a function of exterior wavenumber $k_e$.  The four cubes are aligned as in Figure \ref{fig:cubes}.
In (a) all cubes have low absorption (refractive index $n=1.311 + 2.289 \times 10^{-9}\mathrm{i}$); in (b) all cubes have high absorption (refractive index $n=1.0833 + 0.204 \mathrm{i}$); and in (c) one cube has low absorption, one has high absorption, and the other two have zero absorption (refractive index $n=1.2$).
Other parameters are as in Table \ref{table:iterations_multiple}. 
}
\label{fig:multiple}
\end{figure*}

In Figure \ref{fig:multiple} we compare the performance of the three strong forms $\mathbf{M}^{-1} \mathbf{A}$, $\mathbf{M}^{-1} \mathbf{A}\mathbf{M}^{-1} \mathbf{A}$ and $\mathbf{M}^{-1} \mathbf{D}\mathbf{M}^{-1} \mathbf{A}$, as a function of the exterior wavenumber $k_e$, for (a) low absorption, (b) high absorption, and (c) mixed zero/low/high absorption. 
In all cases the block-diagonally preconditioned form $\mathbf{M}^{-1} \mathbf{D}\mathbf{M}^{-1} \mathbf{A}$ achieved the best performance in terms of overall matvec count, requiring roughly 40\% fewer matvecs than the other two formulations for these four-scatterer configurations.

In Table \ref{table:iterations_multiple2} we compare the performance of the three strong forms as the number of scatterers (cubes) increases. The cubes are aligned as in Figure \ref{fig:cubes}. The GMRES iteration count for $\mathbf{M}^{-1} \mathbf{A} \mathbf{M}^{-1} \mathbf{A}$ and $\mathbf{M}^{-1} \mathbf{D}\mathbf{M}^{-1} \mathbf{A}$ is typically half that for $\mathbf{M}^{-1} \mathbf{A}$. In the light of the complexity calculations \eqref{matvecA}-\eqref{matvecDA}, this suggests that the block-diagonally preconditioned formulation $\mathbf{M}^{-1} \mathbf{D}\mathbf{M}^{-1} \mathbf{A}$ should become more and more efficient compared to the other two methods as the number of scatterers grows, with a theoretical improvement of 50\% fewer matvecs compared to the other two methods in the limit $M\to\infty$. 

\begin{table*}[!t]
\centering
\begin{tabular}{lrrrrrrr}
\toprule
& \multicolumn{3}{c}{$n=1.311 + 2.289 \times 10^{-9}\mathrm{i}$} & &\multicolumn{3}{c}{$n=1.0833 + 0.204 \mathrm{i}$}\\
\cmidrule{2-4} \cmidrule{6-8}
& 4 cubes &  8 cubes & 16 cubes & & 4 cubes &  8 cubes & 16 cubes \\
\midrule
Discrete operator   &   &  & \\
$\mathbf{M}^{-1} \mathbf{A}$ &22 (1840) &24 (7200) &32 (35904)  & &12 (960) &12 (3456)  &13 (14144) \\
$\mathbf{M}^{-1} \mathbf{A}\mathbf{M}^{-1} \mathbf{A}$ &11 (1840) &12 (7200) &16 (35904) & &6 (1040)  &6 (3744) &7 (16320)  \\
$\mathbf{M}^{-1} \mathbf{D}\mathbf{M}^{-1} \mathbf{A}$ &12 (1376) &12 (4288) &16 (19584) & &6 (704) &7 (2528) & 7 (8640) \\
\bottomrule
\end{tabular}
\caption{
Performance of the discrete strong form operators for scattering by arrays of 4, 8 and 16 cubes of side length $0.4$ arranged as in Figure \ref{fig:cubes} at refractive indices $n=1.311 + 2.289 \times 10^{-9}\mathrm{i}$ and $n=1.0833 + 0.204 \mathrm{i}$ and wavenumber $k_e = 20$. Other parameters are as in Table \ref{table:iterations_multiple}. 
}
\label{table:iterations_multiple2}
\end{table*}

\begin{figure*}[t!]
	\centering
    \includegraphics[width = 0.45\textwidth]{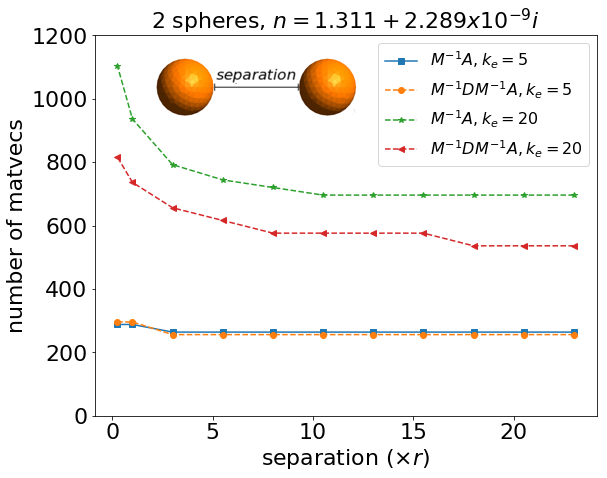}
    \includegraphics[width = 0.45 \textwidth]{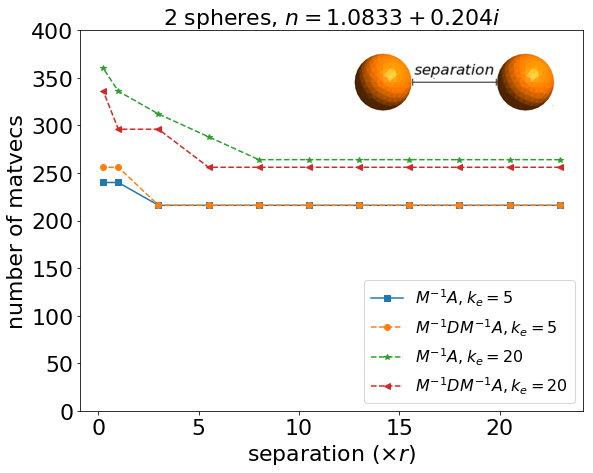}
    \includegraphics[width = 0.45 \textwidth]{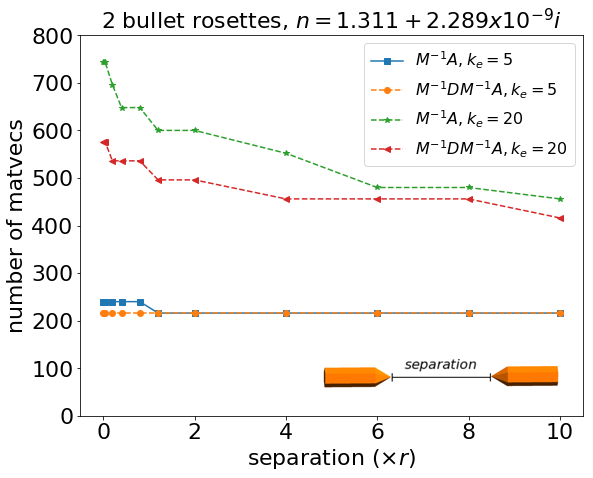}
    \includegraphics[width = 0.45 \textwidth]{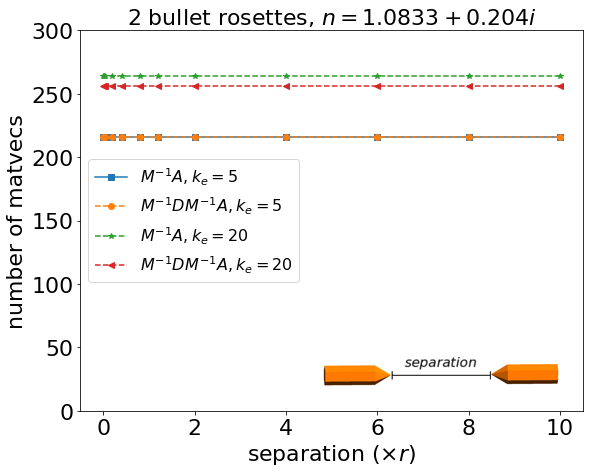}
\caption{Number of matvecs for $\mathbf{M}^{-1} A$ and $\mathbf{M}^{-1} \mathbf{D} \mathbf{M}^{-1} \mathbf{A}$ for scattering by a pair of identical scatterers as a function of their separation. Results are shown for two spheres of radius $0.4$, and for two branches of a bullet rosette. In each case we consider two values of $k_e$ and two refractive indices $n_1=1.0833 + 0.204 \mathrm{i}$ and $n_2=1.311 + 2.289 \times 10^{-9}\mathrm{i}$. The individual size parameters $k_e r$ are $0.4 k_e$ for each sphere and $0.5 k_e$ for each branch of the bullet rosette. The separation is given in terms of the individual particle's radius $r$.
Other parameters are as in Table \ref{table:iterations_multiple}.
}
\label{fig:matvecs_distance}
\end{figure*}

When considering algorithms for scattering by multiple particles it is important to examine the dependence on the distance between the scatterers.
In Figure \ref{fig:matvecs_distance} we investigate how the number of matvecs for the strong form $\mathbf{M}^{-1} \mathbf{A}$ and the block-diagonally preconditioned formulation $\mathbf{M}^{-1} \mathbf{D}\mathbf{M}^{-1} \mathbf{A}$ for a pair of scatterers ($M=2$) depends on the separation between the scatterers. We report results both for spherical scatterers and for two branches of a bullet rosette. 
For the case of low absorption ($n=1.311 + 2.289 \times 10^{-9}\mathrm{i}$), the behaviour depends on the wavenumber $k_e$. For the smaller wavenumber $k_e = 5$ the performance of both operators $\mathbf{M}^{-1} \mathbf{A}$ and $\mathbf{M}^{-1} \mathbf{D}\mathbf{M}^{-1} \mathbf{A}$ is similar. There is no remarkable change in the number of matvecs as the scatterer separation decreases: we observe an increase of 24-48  matvecs (corresponding to 1-2 iterations) for both types of scatterers. For the higher wavenumber $k_e = 20$, $\mathbf{M}^{-1} \mathbf{D}\mathbf{M}^{-1} \mathbf{A}$ performs better than $\mathbf{M}^{-1} \mathbf{A}$ requiring fewer matvecs at any separation level. The number of matvecs required for both types of scatterer is higher, and there is a bigger increase in matvecs as the separation approaches zero.
For the case of high absorption ($n=1.311 + 2.289 \times 10^{-9}\mathrm{i}$), the behaviour depends on both the type of scatterer and the wavenumber $k_e$. For the case of the two branches of a bullet rosette the number of matvecs for each operator is not affected by the separation between them. For the small wavenumber $k_e = 5$, $\mathbf{M}^{-1} \mathbf{A}$ and $\mathbf{M}^{-1} \mathbf{D}\mathbf{M}^{-1} \mathbf{A}$ perform the same but for the higher wavenumber $k_e = 20$, $\mathbf{M}^{-1} \mathbf{D}\mathbf{M}^{-1} \mathbf{A}$ performs better than $\mathbf{M}^{-1} \mathbf{A}$ requiring fewer matvecs. Regarding the two spheres, for the small wavenumber $k_e = 5$, both $\mathbf{M}^{-1} \mathbf{A}$ and $\mathbf{M}^{-1} \mathbf{D}\mathbf{M}^{-1} \mathbf{A}$ perform the same with a modest increase of matvecs as the separation approaches zero. For the higher wavenumber $k_e = 20$, $\mathbf{M}^{-1} \mathbf{D}\mathbf{M}^{-1} \mathbf{A}$ performs better than $\mathbf{M}^{-1} \mathbf{A}$ with both operators requiring an increased number of matvecs as the separation approaches zero. However, even as the separation approaches zero the number of matvecs is not prohibitive. This is in agreement with theoretical results presented in \cite{claeys2012electromagnetic, claeys2012multi}.

These results suggest that for aggregates of scatterers (such as those arising in models of ice crystals in cirrus clouds \cite{baran2009review}), a saving in computational cost might be made by decomposing the aggregate into its constituent parts, treating the aggregate as a multiple scatterer, and applying the block-diagonal preconditioner.  
We investigate this idea further in the next section.

\section{Further numerical examples}\label{sct:ice_crystals}

In this section we present further numerical results, demonstrating that the performance observed for the benchmark problems in the previous section carries over to scattering configurations relevant to light scattering by ice crystals in cirrus clouds \cite{liou2016light,baran2012single,
baran2009review}.

\begin{figure*}[ht!]
\centering
    \begin{subfigure}[t]{\textwidth}
        \hspace{2cm}
        \includegraphics[height = 3cm]{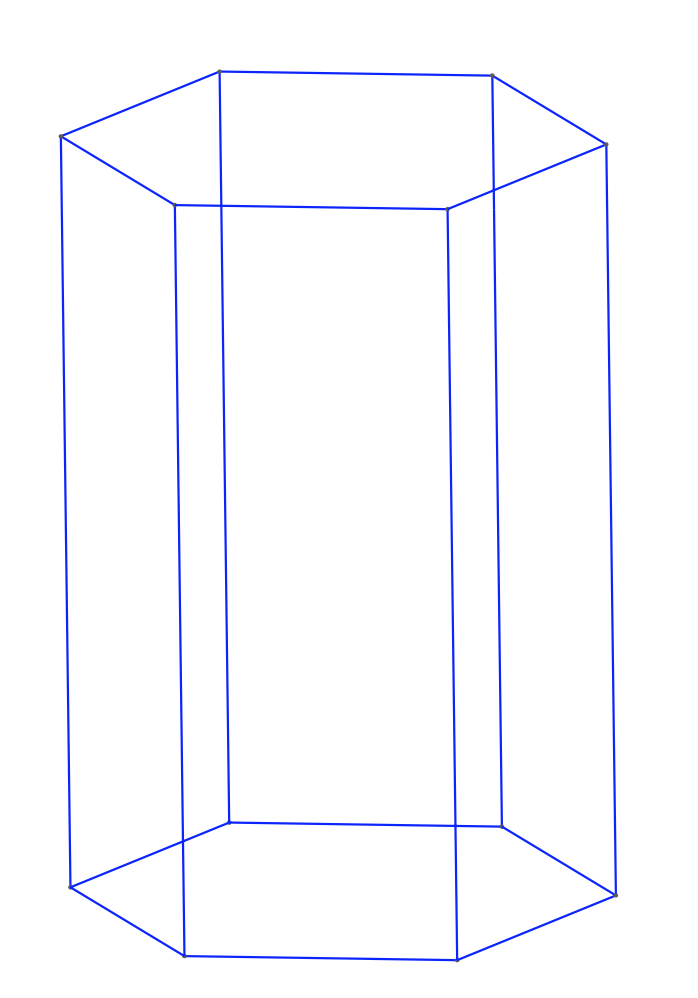} \hspace{4cm} 
        \includegraphics[height = 3.1cm]{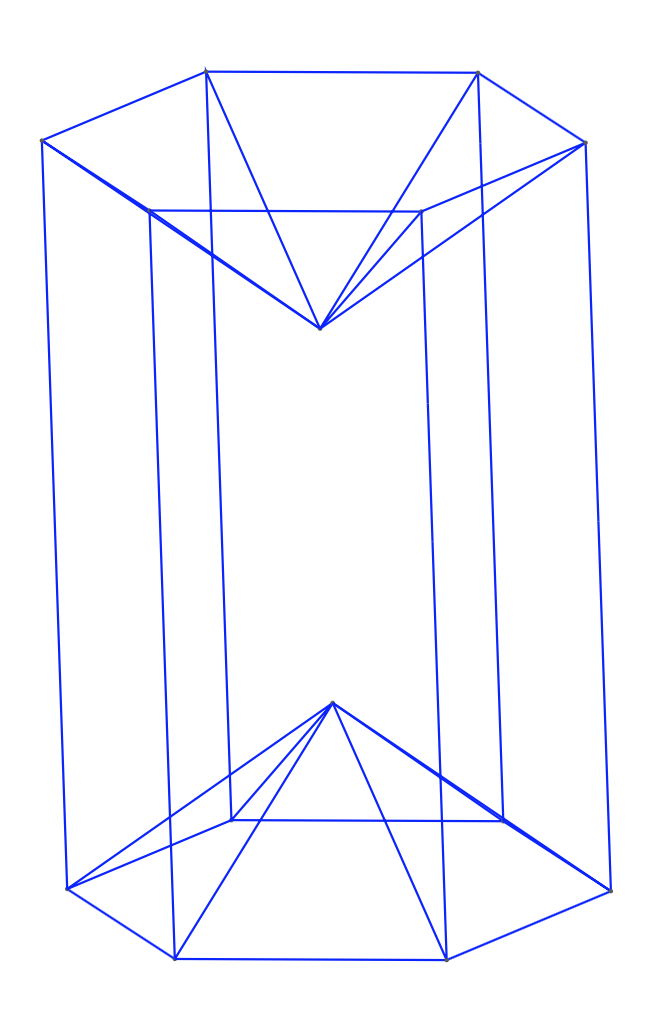} \hspace{3.5cm} 
        \includegraphics[height = 3cm]{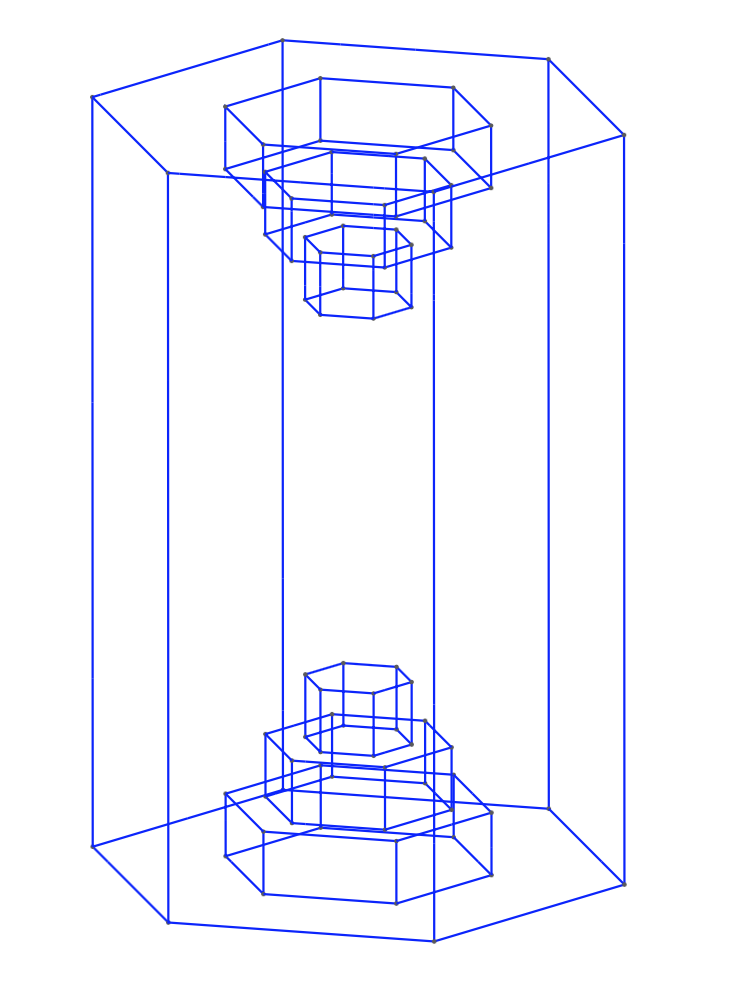}
    \end{subfigure}
    \begin{subfigure}[t]{\textwidth}
        \centering
        \includegraphics[width = 0.32 \textwidth]{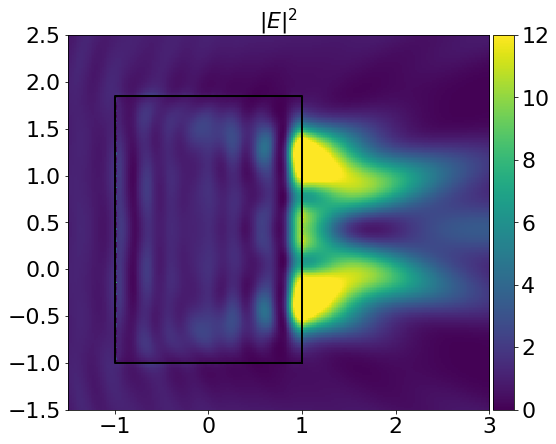} 
        \includegraphics[width = 0.32 \textwidth]{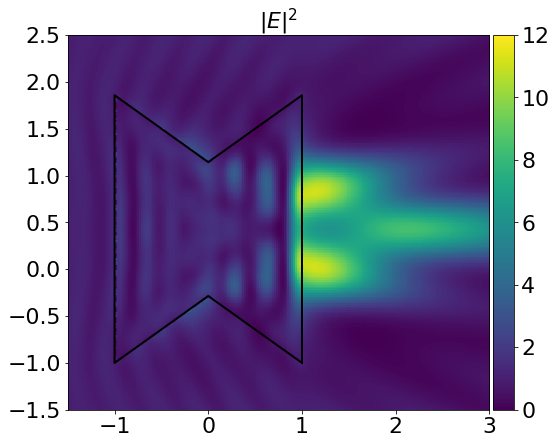} 
        \includegraphics[width = 0.32 \textwidth]{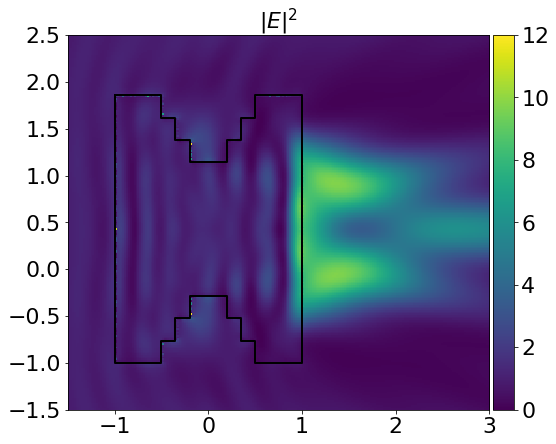}
    \end{subfigure}
    \caption{Squared magnitude $|\mathbf{E}|^2$ of the electric field for single-particle scattering $(M=1)$ by a hexagonal column of increased complexity (with no cavities, with a conventional cavity and with a stepped-cavity) restricted to the $xz$ plane. The incident wave is $\mathbf{E}^{inc}(\mathbf{x})=\mathbf{p} \mathrm{e}^{\mathrm{i} k_e \mathbf{d} \cdot \mathbf{x}}$, with $\mathbf{d} = (1,0,0)^T$ and $\mathbf{p} = (0,0,1)^T$, with $k_e = 10$, $n=1.311 + 2.289 \times 10^{-9}\mathrm{i}$, $k_1 = n k_e$, $\mu_e = \mu_1 = 1$. The size parameter $k_e r$ is 30. The mesh size is $h = 2\pi/(10 k_e)$.}
    \label{fig:hex_col_results}
\end{figure*}

We first present results for single-particle scattering problems for three hexagonal columns of increasing geometric complexity: without cavities, with ``conventional'' cavities, and with stepped cavities \cite{smith2015cloud} (see Figure \ref{fig:hex_col_results}). The applicability of Bempp to these particles was already demonstrated in \cite{groth2015boundary}, where far-field scattering properties for the three cases were compared and, in the case of standard hexagonal columns, validated against corresponding T-matrix calculations. In Figure \ref{fig:hex_col_results} we show plots of the squared magnitude $|\mathbf{E}|^2$ of the electric field inside and outside the particles, restricted to a plane parallel to the column axes, for a particular scattering configuration. 
These plots supplement those already presented in \cite{groth2015boundary}, demonstrating that Bempp can easily generate near field as well as far field plots. 
The near field interference patterns for the three particles are similar with a stronger focusing for the hexagonal column without cavities. Between the three scatterers, the column with the stepped cavities exhibits the weakest focusing.

In Table \ref{table:hex_columns} we report the corresponding GMRES iteration and matvec counts for the two discrete strong forms $\mathbf{M}^{-1} \mathbf{A}$ and $\mathbf{M}^{-1} \mathbf{A}\mathbf{M}^{-1} \mathbf{A}$. 
Both operators perform the same in terms of overall matvec count, despite the mass-matrix preconditioner needing twice the number of GMRES iterations.
Interestingly, the number of matvecs and iterations required is roughly the same for all three examples, which demonstrates the effectiveness of BEM for complex particle geometries. 

\begin{table}[t!]
\centering
\begin{tabular}{lrrr}
\toprule
& \multicolumn{3}{c}{$n=1.311 + 2.289 \times 10^{-9}\mathrm{i}$} \\
\cmidrule{2-4} 
& hexagonal &  with & with \\
& column &  conventional & stepped \\
&  &  cavity & cavity \\
\midrule
Discrete operator   &   &  & \\
$\mathbf{M}^{-1} \mathbf{A}$ & 34 (280) & 30 (248) & 34 (280)\\
$\mathbf{M}^{-1} \mathbf{A}\mathbf{M}^{-1} \mathbf{A}$ &17 (280) &15 (248) &17 (280)\\
\bottomrule
\end{tabular}
\caption{Number of GMRES iterations and total matvec count (in brackets) for the different discrete formulations for single-particle scattering  ($M=1$) by the particles of Figure \ref{fig:hex_col_results}. The parameters are the same as in Figure \ref{fig:hex_col_results}.}
\label{table:hex_columns}
\end{table}

\begin{figure*}[t!]
\centering
    \begin{subfigure}[t]{\textwidth}
        \centering
        \includegraphics[height = 4cm]{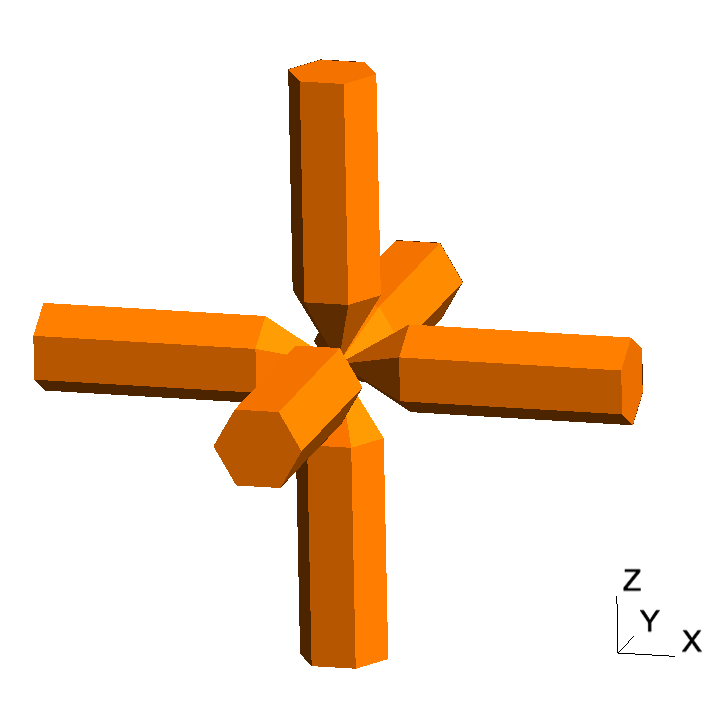} \hspace{2.5cm}
        \includegraphics[height = 4cm]{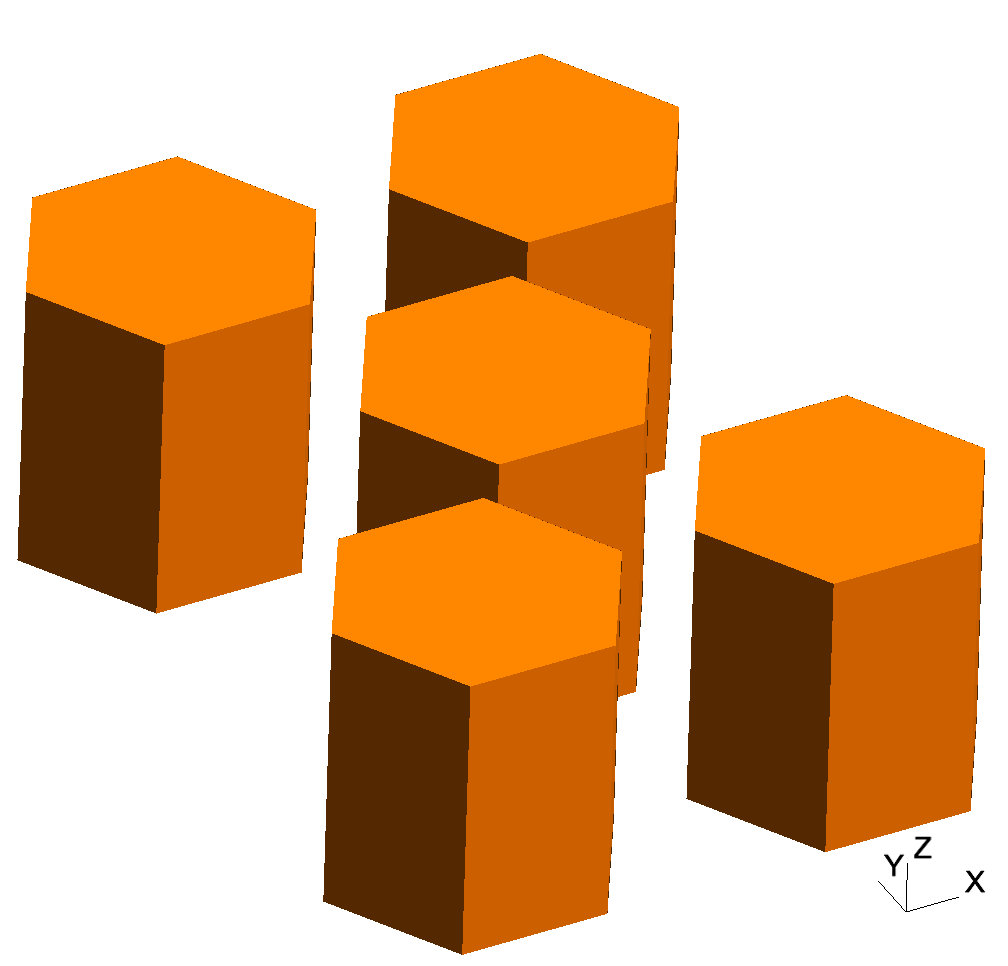} \hspace{3cm}
        \includegraphics[height = 4cm]{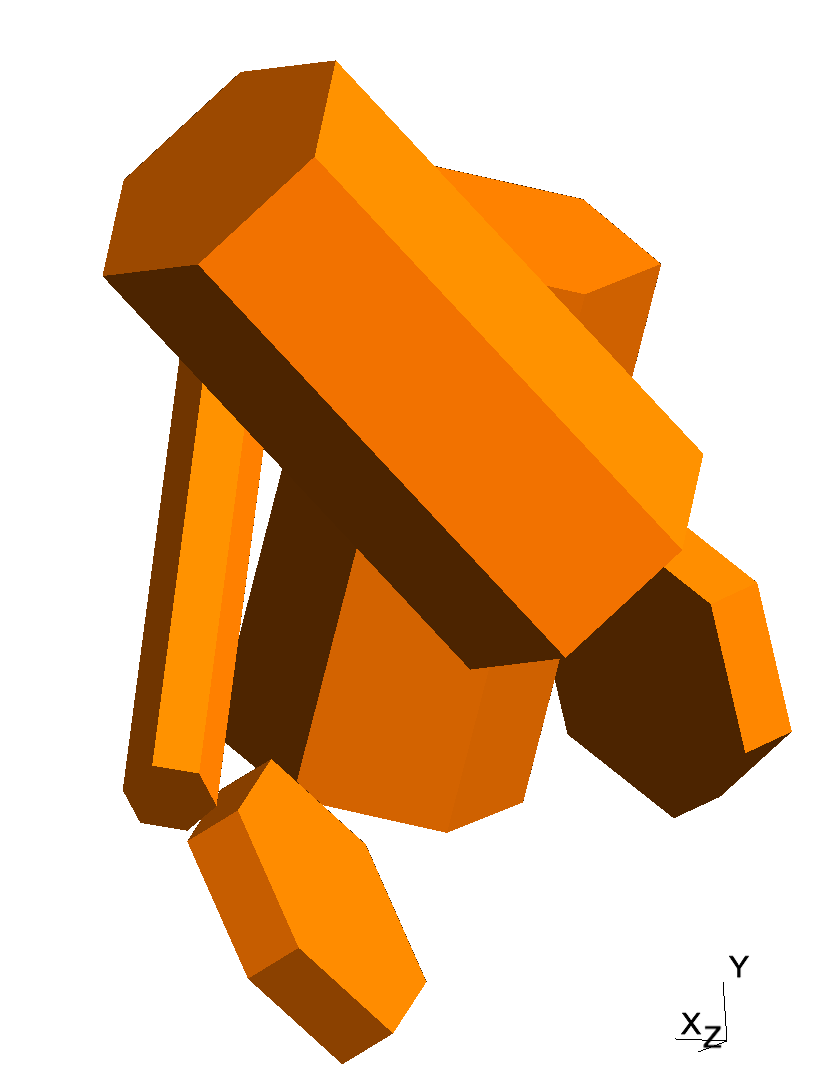}
    \end{subfigure}
    \begin{subfigure}[t]{\textwidth}
        \centering
        \includegraphics[width = 0.32 \textwidth]{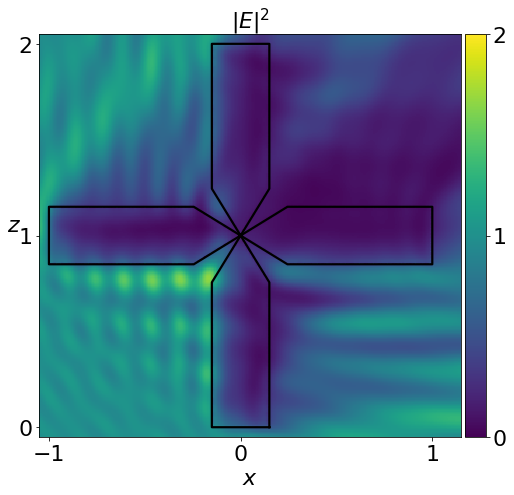}
        \includegraphics[width = 0.33 \textwidth]{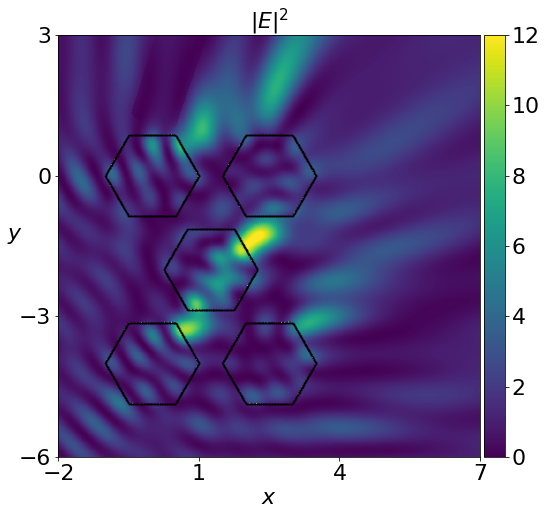}
        \includegraphics[width = 0.33 \textwidth]{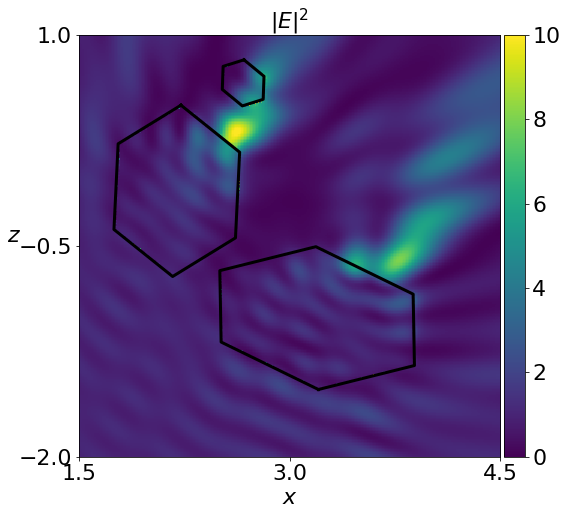}
    \end{subfigure}
    \caption{Squared magnitude $|\mathbf{E}|^2$ of the electric field for scattering by multiple particles $(M>1)$ restricted to different planes. 
    In the case of the bullet-rosette with 6 branches ($M=6$) the incident wave is $\mathbf{E}^{inc}(\mathbf{x})=\mathbf{p} \mathrm{e}^{\mathrm{i} k_e \mathbf{d} \cdot \mathbf{x}}$, with $\mathbf{d} = (\sqrt{3}/2,0,1/2)^T$ and $\mathbf{p} = (0,1,0)^T$, with $k_e = 25$, $n = 1.0833 + 0.204 \mathrm{i}$, $k_m = n k_e$, and $\mu_e = \mu_m = 1$, for $m = 1, \ldots, 6$. The size parameter of the aggregate is $k_e r = 25 \sqrt{2}$. 
    For the array of five hexagonal columns ($M=5$) the incident wave is $\mathbf{E}^{inc}(\mathbf{x})=\mathbf{p} \mathrm{e}^{\mathrm{i} k_e \mathbf{d} \cdot \mathbf{x}}$, with $\mathbf{d} = (1/\sqrt{2},1/\sqrt{2},0)$ and $\mathbf{p} = (0,0,1)^T$, with $k_e = 5$, $n = 1.311 + 2.289 \times 10^{-9}\mathrm{i}$, $k_m = n k_e$, and $\mu_e = \mu_m = 1$, for $m = 1, \ldots, 5$. The size parameter of the aggregate is $k_e r = 14$. For the aggregate of randomly oriented hexagonal columns and plates ($M=5$) the incident wave is $\mathbf{E}^{inc}(\mathbf{x})=\mathbf{p} \mathrm{e}^{\mathrm{i} k_e \mathbf{d} \cdot \mathbf{x}}$, with $\mathbf{d} = (1/\sqrt{2},0,1/\sqrt{2})$ and $\mathbf{p} = (0,1,0)^T$, with $k_e = 15$, $n = 1.311 + 2.289 \times 10^{-9}\mathrm{i}$, $k_m = n k_e$, and $\mu_e = \mu_m = 1$, for $m = 1, \ldots, 5$. The size parameter of the aggregate is $k_e r = 49$. The mesh size in all cases is $h = 2\pi/(10 k_e)$.}
    \label{fig:aggregates_results}
\end{figure*}

We now present results for three problems of scattering by multiple particles involving ensembles or aggregates of simpler particles: a six-branch bullet rosette, an aligned array of five hexagonal columns of the same dimensions, and a random aggregate of five hexagonal columns/plates of differing dimensions (see Figure \ref{fig:aggregates_results}). In Figure \ref{fig:aggregates_results} we show the scatterer geometries, along with the 
fields obtained using the BEM formulation for particular scattering configurations (with parameters detailed in the figure caption). For the bullet rosette we present the near field plot of $|\mathbf{E}|^2$ restricted to the (vertical) $xz$-plane, for the aligned array the near field plot of $|\mathbf{E}|^2$ restricted to the (horizontal) $xy$-plane, and for the random aggregate the near field plot of $|\mathbf{E}|^2$ restricted to the (vertical) $xz$-plane.
The calculation was carried out on a high-spec desktop machine with assembly taking roughly 23 minutes and the linear solve taking roughly 7 minutes.

The corresponding GMRES iteration and matvec counts are reported in Table \ref{table:aggregates}. We note that for all three multi-particle configurations, the block-diagonally Calder\'on-preconditioned approach $\mathbf{M}^{-1} \mathbf{D}\mathbf{M}^{-1} \mathbf{A}$ provides a significant saving in overall matvec count compared to the other two approaches, even though the constituent particles in the bullet rosette and random aggregate are close enough to touch each other at certain points. 

\begin{table}[t!]
\centering
\begin{tabular}{lrrr}
\toprule
& 6-branch &  5 hex. & random \\
& bullet &  columns & columns \\
& rosette & & \\
\midrule
Discrete operator   &   &  & \\
$\mathbf{M}^{-1} \mathbf{A}$ & 15 (2520) & 52 (6480)  & 38 (4680)\\
$\mathbf{M}^{-1} \mathbf{A}\mathbf{M}^{-1} \mathbf{A}$ & 8 (2856) & 25 (6360) & 19 (4680)\\
$\mathbf{M}^{-1} \mathbf{D}\mathbf{M}^{-1} \mathbf{A}$ & 8 (1776) & 23 (3880) & 19 (3080)\\
\bottomrule
\end{tabular}
\caption{Number of GMRES iterations and total matvec count (in brackets) for the different discrete formulations for scattering  by the particles ($M>1$) of Figure \ref{fig:aggregates_results}. The parameters are the same as in Figure \ref{fig:aggregates_results}.}
\label{table:aggregates}
\end{table}

\section{Conclusion}\label{sct:conclusion}
We have carried out a detailed study into the performance of various operator-based preconditioning strategies for Galerkin BEM discretizations of the PMCHWT boundary integral equation for electromagnetic scattering by absorbing dielectric particles of different shapes, sizes and refractive indices. Specifically, we considered the weak and strong discrete forms of (i) the PMCHWT operator $\bm{\mathcal{A}}$, (ii) its Calder\'on-preconditioned square $\bm{\mathcal{A}}^2$, and, for multiple scatterers, (iii) a block-diagonally Calder\'on-preconditioned operator $\bm{\mathcal{D}}\bm{\mathcal{A}}$,  in which only the self-interaction blocks are preconditioned. Numerical performance was measured in the number of matrix-vector products (``matvecs'') incurred when the corresponding linear system was solved using GMRES. 

Overall we found that the strong forms (which involve multiplications by inverse mass-matrices) required significantly fewer matvecs than their weak counterparts. The strong form of the Calder\'on-preconditioned operator $\bm{\mathcal{A}}^2$ in general led to a reduction by 50\% in the number of GMRES iterations required compared to the strong form of $\bm{\mathcal{A}}$ (i.e., simple mass-matrix preconditioning). But in terms of total matvecs the increased cost per iteration outweighed this gain. Hence for single scattering applications we found that simple mass-matrix preconditioning was more effective than Calder\'on preconditioning.

For problems of scattering by multiple particles  we found that a saving of up to 50\% in total matvec count compared to the mass-matrix preconditioner could be achieved by using the strong form of the block-diagonal Calder\'on preconditioned operator $\bm{\mathcal{D}}\bm{\mathcal{A}}$. Provided that the particles were sufficiently absorbing, this gain in performance was found to hold even when the particles were close, or even touching each other. 

Using the BEM software library Bempp the preconditioners were applied to various scattering configurations relevant to the scattering of light by atmospheric ice crystals including hexagonal columns with conventional and stepped cavities, bullet rosettes and a random aggregate of hexagonal columns, demonstrating the applicability of the methods for practical scattering simulations. 

\section*{Acknowledgements}
The work of the first author was supported by the Natural Environment Research Council and the UK Met Office (CASE PhD studentship to A.\ Kleanthous, grant NE/N008111/1).

\section*{References}
\bibliography{mybibfile}

\begin{thebibliography}{10}
\expandafter\ifx\csname url\endcsname\relax
  \def\url#1{\texttt{#1}}\fi
\expandafter\ifx\csname urlprefix\endcsname\relax\def\urlprefix{URL }\fi
\expandafter\ifx\csname href\endcsname\relax
  \def\href#1#2{#2} \def\path#1{#1}\fi

\bibitem{liou2016light}
K.-N. Liou, P.~Yang, {Light scattering by ice crystals: fundamentals and
  applications}, Cambridge University Press, 2016.

\bibitem{baran2012single}
A.~J. Baran, {From the single-scattering properties of ice crystals to climate
  prediction: A way forward}, Atmos. Res. 112 (2012) 45--69.

\bibitem{baran2009review}
A.~J. Baran, {A review of the light scattering properties of cirrus}, J. Quant.
  Spectrosc. Ra. 110~(14) (2009) 1239--1260.

\bibitem{mishchenko2014electromagnetic}
M.~I. Mishchenko, {Electromagnetic scattering by particles and particle groups:
  an introduction}, Cambridge University Press, 2014.

\bibitem{draine1994discrete}
B.~T. Draine, P.~J. Flatau, {Discrete-dipole approximation for scattering
  calculations}, JOSA A 11~(4) (1994) 1491--1499.

\bibitem{yurkin2007discrete}
M.~A. Yurkin, A.~G. Hoekstra, {The discrete dipole approximation: an overview
  and recent developments}, J. Quant. Spectrosc. Ra. 106~(1-3) (2007) 558--589.

\bibitem{yang1996finite}
P.~Yang, K.~Liou, {Finite-difference time domain method for light scattering by
  small ice crystals in three-dimensional space}, JOSA A 13~(10) (1996)
  2072--2085.

\bibitem{yang2000finite}
P.~Yang, K.~Liou, {Finite difference time domain method for light scattering by
  nonspherical and inhomogeneous particles}, in: Light Scattering by
  Nonspherical Particles: Theory, Measurements, and Applications, Vol.~1,
  Academic, 2000, p. 174.

\bibitem{yurkin2007systematic}
M.~A. Yurkin, A.~G. Hoekstra, R.~S. Brock, J.~Q. Lu, {Systematic comparison of
  the discrete dipole approximation and the finite difference time domain
  method for large dielectric scatterers}, Opt. Express 15~(26) (2007)
  17902--17911.

\bibitem{liu2012comparison}
C.~Liu, L.~Bi, R.~L. Panetta, P.~Yang, M.~A. Yurkin, {Comparison between the
  pseudo-spectral time domain method and the discrete dipole approximation for
  light scattering simulations}, Opt. Express 20~(15) (2012) 16763--16776.

\bibitem{havemann2001extension}
S.~Havemann, A.~Baran, {Extension of T-matrix to scattering of electromagnetic
  plane waves by non-axisymmetric dielectric particles: application to
  hexagonal ice cylinders}, J. Quant. Spectrosc. Ra. 70~(2) (2001) 139--158.

\bibitem{kahnert2013t}
M.~Kahnert, {The T-matrix code Tsym for homogeneous dielectric particles with
  finite symmetries}, J. Quant. Spectrosc. Ra. 123 (2013) 62--78.

\bibitem{baran2001calculation}
A.~J. Baran, P.~Yang, S.~Havemann, {Calculation of the single-scattering
  properties of randomly oriented hexagonal ice columns: a comparison of the
  T-matrix and the finite-difference time-domain methods}, Appl. Optics 40~(24)
  (2001) 4376--4386.

\bibitem{bi2013efficient}
L.~Bi, P.~Yang, G.~W. Kattawar, M.~I. Mishchenko, {Efficient implementation of
  the invariant imbedding T-matrix method and the separation of variables
  method applied to large nonspherical inhomogeneous particles}, J. Quant.
  Spectrosc. Ra. 116 (2013) 169--183.

\bibitem{bi2014accurate}
L.~Bi, P.~Yang, {Accurate simulation of the optical properties of atmospheric
  ice crystals with the invariant imbedding T-matrix method}, J. Quant.
  Spectrosc. Ra. 138 (2014) 17--35.

\bibitem{muinonen1989scattering}
K.~Muinonen, {Scattering of light by crystals: a modified Kirchhoff
  approximation}, Appl. Optics 28~(15) (1989) 3044--3050.

\bibitem{muinonen1996light}
K.~Muinonen, T.~Nousiainen, P.~Fast, K.~Lumme, J.~Peltoniemi, {Light scattering
  by Gaussian random particles: ray optics approximation}, J. Quant. Spectrosc.
  Ra. 55~(5) (1996) 577--601.

\bibitem{macke1996single}
A.~Macke, J.~Mueller, E.~Raschke, {Single scattering properties of atmospheric
  ice crystals}, J. Atmos. Sci. 53~(19) (1996) 2813--2825.

\bibitem{yang1996geometric}
P.~Yang, K.~Liou, {Geometric-optics--integral-equation method for light
  scattering by nonspherical ice crystals}, Appl. Optics 35~(33) (1996)
  6568--6584.

\bibitem{mishchenko1998incorporation}
M.~I. Mishchenko, A.~Macke, {Incorporation of physical optics effects and
  computation of the Legendre expansion for ray-tracing phase functions
  involving $\delta$-function transmission}, J. Geophys. Res. - Atmos. 103~(D2)
  (1998) 1799--1805.

\bibitem{liou2002introduction}
K.-N. Liou, {An introduction to atmospheric radiation}, Elsevier, 2002.

\bibitem{borovoi2003scattering}
A.~G. Borovoi, I.~A. Grishin, Scattering matrices for large ice crystal
  particles, JOSA A 20~(11) (2003) 2071--2080.

\bibitem{hesse2008modelling}
E.~Hesse, Modelling diffraction during ray tracing using the concept of energy
  flow lines, J. Quant. Spectrosc. Ra. 109~(8) (2008) 1374--1383.

\bibitem{bi2009simulation}
L.~Bi, P.~Yang, G.~W. Kattawar, B.~A. Baum, Y.~X. Hu, D.~M. Winker, R.~S.
  Brock, J.~Q. Lu, {Simulation of the color ratio associated with the
  backscattering of radiation by ice particles at the wavelengths of 0.532 and
  1.064 $\mu$m}, J. Geophys. Res. - Atmos. 114~(D4).

\bibitem{bi2011scattering}
L.~Bi, P.~Yang, G.~W. Kattawar, Y.~Hu, B.~A. Baum, {Scattering and absorption
  of light by ice particles: solution by a new physical-geometric optics hybrid
  method}, J. Quant. Spectrosc. Ra. 112~(9) (2011) 1492--1508.

\bibitem{hesse2012modelling}
E.~Hesse, A.~Macke, S.~Havemann, A.~Baran, Z.~Ulanowski, P.~H. Kaye, {Modelling
  diffraction by facetted particles}, J. Quant. Spectrosc. Ra. 113~(5) (2012)
  342--347.

\bibitem{bi2014assessment}
L.~Bi, P.~Yang, C.~Liu, B.~Yi, B.~A. Baum, B.~Van~Diedenhoven, H.~Iwabuchi,
  {Assessment of the accuracy of the conventional ray-tracing technique:
  Implications in remote sensing and radiative transfer involving ice clouds},
  J. Quant. Spectrosc. Ra. 146 (2014) 158--174.

\bibitem{martin2006multiple}
P.~A. Martin, {Multiple scattering: interaction of time-harmonic waves with N
  obstacles}, no. 107, Cambridge University Press, 2006.

\bibitem{koc1998calculation}
S.~Koc, W.~Chew, {Calculation of acoustical scattering from a cluster of
  scatterers}, The Journal of the Acoustical Society of America 103~(2) (1998)
  721--734.

\bibitem{gumerov2005computation}
N.~A. Gumerov, R.~Duraiswami, {Computation of scattering from clusters of
  spheres using the fast multipole method}, The Journal of the Acoustical
  Society of America 117~(4) (2005) 1744--1761.

\bibitem{zhang2007fast}
Y.~J. Zhang, E.~P. Li, {Fast multipole accelerated scattering matrix method for
  multiple scattering of a large number of cylinders}, Progress In
  Electromagnetics Research 72 (2007) 105--126.

\bibitem{bremer2015high}
J.~Bremer, A.~Gillman, P.-G. Martinsson, {A high-order accurate accelerated
  direct solver for acoustic scattering from surfaces}, BIT Numerical
  Mathematics 55~(2) (2015) 367--397.

\bibitem{ganesh2015efficient}
M.~Ganesh, S.~Hawkins, {An efficient $\mathcal{O}(N)$ algorithm for computing
  $\mathcal{O}(N^2)$ acoustic wave interactions in large $N$-obstacle three
  dimensional configurations}, BIT Numerical Mathematics 55~(1) (2015)
  117--139.

\bibitem{groth2015boundary}
S.~Groth, A.~Baran, T.~Betcke, S.~Havemann, W.~{\'S}migaj, {The boundary
  element method for light scattering by ice crystals and its implementation in
  BEM++}, J. Quant. Spectrosc. Ra. 167 (2015) 40--52.

\bibitem{baran2017application}
A.~J. Baran, S.~P. Groth, {The application of the boundary element method in
  BEM++ to small extreme Chebyshev ice particles and the remote detection of
  the ice crystal number concentration of small atmospheric ice particles}, J.
  Quant. Spectrosc. Ra. 198 (2017) 68--80.

\bibitem{yu2017electromagnetic}
M.~P. Yu, Y.~P. Han, Z.~W. Cui, A.~T. Chen, {Electromagnetic scattering by
  multiple dielectric particles under the illumination of unpolarized
  high-order Bessel vortex beam}, J. Quant. Spectrosc. Ra. 195 (2017) 107--113.

\bibitem{mano2000exact}
Y.~Mano, {Exact solution of electromagnetic scattering by a three-dimensional
  hexagonal ice column obtained with the boundary-element method}, Appl. Optics
  39~(30) (2000) 5541--5546.

\bibitem{nakajima2009development}
T.~Y. Nakajima, T.~Nakajima, K.~Yoshimori, S.~K. Mishra, S.~N. Tripathi,
  {Development of a light scattering solver applicable to particles of
  arbitrary shape on the basis of the surface-integral equations method of
  M{\"u}ller type. I. Methodology, accuracy of calculation, and electromagnetic
  current on the particle surface}, Appl. Optics 48~(19) (2009) 3526--3536.

\bibitem{smigaj2015solving}
W.~{\'S}migaj, T.~Betcke, S.~Arridge, J.~Phillips, M.~Schweiger, {Solving
  boundary integral problems with BEM++}, ACM T. Math. Software (TOMS) 41~(2)
  (2015) 6.

\bibitem{groth2018hybrid}
S.~P. Groth, D.~P. Hewett, S.~Langdon, A hybrid numerical--asymptotic boundary
  element method for high frequency scattering by penetrable convex polygons,
  Wave Motion 78 (2018) 32--53.

\bibitem{poggio1970integral}
A.~J. Poggio, E.~K. Miller, {Integral equation solutions of three-dimensional
  scattering problems}, MB Assoc., 1970.

\bibitem{wu1977scattering}
T.-K. Wu, L.~L. Tsai, {Scattering from arbitrarily-shaped lossy dielectric
  bodies of revolution}, Radio Sci. 12~(5) (1977) 709--718.

\bibitem{mautz1977electromagnetic}
J.~R. Mautz, R.~F. Harrington, {Electromagnetic scattering from a homogeneous
  body of revolution}, Tech. rep., Syracuse Univ. NY Dept. of Electrical and
  Computer Engineering (1977).

\bibitem{harrington1989boundary}
R.~F. Harrington, {Boundary integral formulations for homogeneous material
  bodies}, J. Electromagnet. Wave. 3~(1) (1989) 1--15.

\bibitem{saad1986gmres}
Y.~Saad, M.~H. Schultz, {GMRES: A generalized minimal residual algorithm for
  solving nonsymmetric linear systems}, SIAM J. Sci. Stat. Comp. 7~(3) (1986)
  856--869.

\bibitem{contopanagos2002well}
H.~Contopanagos, B.~Dembart, M.~Epton, J.~J. Ottusch, V.~Rokhlin, J.~L. Visher,
  S.~M. Wandzura, {Well-conditioned boundary integral equations for
  three-dimensional electromagnetic scattering}, IEEE T. Antenn. Propag.
  50~(12) (2002) 1824--1830.

\bibitem{christiansen2002preconditioner}
S.~H. Christiansen, J.-C. N{\'e}d{\'e}lec, {A preconditioner for the electric
  field integral equation based on Calderon formulas}, SIAM journal on
  numerical analysis 40~(3) (2002) 1100--1135.

\bibitem{antoine2008integral}
X.~Antoine, Y.~Boubendir, {An integral preconditioner for solving the
  two-dimensional scattering transmission problem using integral equations},
  International Journal of Computer Mathematics 85~(10) (2008) 1473--1490.

\bibitem{bagci2009calderon}
H.~Bagci, F.~P. Andriulli, K.~Cools, F.~Olyslager, E.~Michielssen, {A
  Calder{\'o}n multiplicative preconditioner for the combined field integral
  equation}, IEEE Transactions on Antennas and Propagation 57~(10) (2009)
  3387--3392.

\bibitem{cools2011calderon}
K.~Cools, F.~P. Andriulli, E.~Michielssen, {A Calder{\'o}n multiplicative
  preconditioner for the PMCHWT integral equation}, IEEE T. Antenn. Propag.
  59~(12) (2011) 4579--4587.

\bibitem{yan2010comparative}
S.~Yan, J.-M. Jin, Z.~Nie, {A comparative study of Calder{\'o}n preconditioners
  for PMCHWT equations}, IEEE T. Antenn. Propag. 58~(7) (2010) 2375--2383.

\bibitem{niino2012calderon}
K.~Niino, N.~Nishimura, {Calder{\'o}n preconditioning approaches for PMCHWT
  formulations for Maxwell's equations}, Int. J. Numer. Model. El. 25~(5-6)
  (2012) 558--572.

\bibitem{boubendir2015integral}
Y.~Boubendir, O.~Bruno, D.~Levadoux, C.~Turc, {Integral equations requiring
  small numbers of Krylov-subspace iterations for two-dimensional smooth
  penetrable scattering problems}, Applied Numerical Mathematics 95 (2015)
  82--98.

\bibitem{kirby2010functional}
R.~C. Kirby, From functional analysis to iterative methods, SIAM Review 52~(2)
  (2010) 269--293.

\bibitem{betcke2017product}
T.~Betcke, M.~Scroggs, W.~Smigaj, {Product algebras for Galerkin
  discretisations of boundary integral operators and their applications}, arXiv
  preprint arXiv:1711.10607.

\bibitem{buffa2003galerkin}
A.~Buffa, R.~Hiptmair, {Galerkin boundary element methods for electromagnetic
  scattering}, in: Topics in computational wave propagation, Springer, 2003,
  pp. 83--124.

\bibitem{kirsch2015mathematical}
A.~Kirsch, F.~Hettlich, {The Mathematical Theory of Time-Harmonic Maxwell's
  Equations}, Applied Mathematical Sciences 190.

\bibitem{muller2013foundations}
C.~M{\"u}ller, {Foundations of the mathematical theory of electromagnetic
  waves}, Vol. 155, Springer Science \& Business Media, 2013.

\bibitem{nedelec2001acoustic}
J.-C. N{\'e}d{\'e}lec, {Acoustic and electromagnetic equations: integral
  representations for harmonic problems}, Vol. 144, Springer Science \&
  Business Media, 2001.

\bibitem{rwg1982}
S.~Rao, D.~Wilton, A.~Glisson, {Electromagnetic scattering by surfaces of
  arbitrary shape}, IEEE T. Antenn. Propag. 30~(3) (1982) 409--418.

\bibitem{buffa2007dual}
A.~Buffa, S.~Christiansen, {A dual finite element complex on the barycentric
  refinement}, Math. Comp. 76~(260) (2007) 1743--1769.

\bibitem{scroggs2017software}
M.~W. Scroggs, T.~Betcke, E.~Burman, W.~{\'S}migaj, E.~van't Wout, {Software
  frameworks for integral equations in electromagnetic scattering based on
  Calder{\'o}n identities}, Comput. Math. Appl. 74~(11) (2017) 2897--2914.

\bibitem{hackbusch2015hierarchical}
W.~Hackbusch, Hierarchical matrices: algorithms and analysis, Vol.~49,
  Springer, 2015.

\bibitem{gumerov2005fast}
N.~A. Gumerov, R.~Duraiswami, Fast multipole methods for the Helmholtz equation
  in three dimensions, Elsevier, 2005.

\bibitem{engquist2018approximate}
B.~Engquist, H.~Zhao, Approximate separability of the green's function of the
  helmholtz equation in the high frequency limit, Comm. Pure Appl. Math.

\bibitem{betcke2017computationally}
T.~Betcke, E.~van't Wout, P.~G{\'e}lat, Computationally efficient boundary
  element methods for high-frequency helmholtz problems in unbounded domains,
  in: Modern Solvers for Helmholtz Problems, Springer, 2017, pp. 215--243.

\bibitem{warren2008optical}
S.~G. Warren, R.~E. Brandt, {Optical constants of ice from the ultraviolet to
  the microwave: A revised compilation}, J. Geophys. Res. - Atmos. 113~(D14).

\bibitem{claeys2012electromagnetic}
X.~Claeys, R.~Hiptmair, {Electromagnetic scattering at composite objects: a
  novel multi-trace boundary integral formulation}, ESAIM: Mathematical
  Modelling and Numerical Analysis 46~(6) (2012) 1421--1445.

\bibitem{claeys2012multi}
X.~Claeys, R.~Hiptmair, C.~Jerez-Hanckes, {Multi-trace boundary integral
  equations}, Direct and inverse problems in wave propagation and applications
  14 (2012) 51--100.

\bibitem{smith2015cloud}
H.~R. Smith, P.~J. Connolly, A.~J. Baran, E.~Hesse, A.~R. Smedley, A.~R. Webb,
  {Cloud chamber laboratory investigations into scattering properties of hollow
  ice particles}, J. Quant. Spectrosc. Ra. 157 (2015) 106--118.

\end{thebibliography}

\end{document}